\documentclass[10pt,journal,twocolumn]{IEEEtran}

\usepackage{ifpdf}
\usepackage{steinmetz}
\usepackage{amsmath}
\usepackage{amssymb}    
\usepackage[inline]{enumitem}
\usepackage{threeparttable}
\newcommand\MYhyperrefoptions{bookmarks=true,bookmarksnumbered=true,
pdfpagemode={UseOutlines},plainpages=false,pdfpagelabels=true,
colorlinks=true,linkcolor={black},citecolor={black},urlcolor={black},
pdftitle={threephasemodeling},%<!CHANGE!
pdfsubject={},%<!CHANGE!
pdfauthor={Mohammadhafez Bazrafshan and Nikolaos Gatsis},%<!CHANGE!
pdfkeywords={}}%<^!CHANGE!

\ifCLASSINFOpdf
\usepackage[\MYhyperrefoptions,pdftex, draft]{hyperref}
\else
\usepackage[\MYhyperrefoptions,breaklinks=true,dvips]{hyperref}
\usepackage{breakurl}
\fi
\newtheorem{theorem}{Theorem}
\newtheorem{remark}{Remark}

\newtheorem{lemma}{Lemma}
\newtheorem{corollary}{Corollary}

\newcommand{\mr}[1]{\mathrm{#1}}
\newcommand{\mc}[1]{\mathcal{#1}}
\newcommand{\mb}[1]{\mathbf{#1}} % general symbols
\newcommand{\mbb}[1]{\mathbb{#1}}
\newcommand{\bmat}[1]{\begin{bmatrix} #1 \end{bmatrix}}
\newcommand{\smat}[1]{\left[\begin{smallmatrix} #1 \end{smallmatrix}\right]}

\newcommand{\eq}[2]{
	\begin{IEEEeqnarray}{#1}
		#2
	\end{IEEEeqnarray}}

\usepackage{multirow}

\ifCLASSOPTIONcompsoc
\usepackage[caption=false,font=normalsize,labelfon
t=sf,textfont=sf]{subfig}
\else
\usepackage[caption=false,font=footnotesize]{subfi
g}
\fi

\ifCLASSINFOpdf
\usepackage[pdftex]{graphicx}
 \usepackage{epstopdf}
\DeclareGraphicsExtensions{.pdf,.eps}
\else
 \usepackage[dvips]{graphicx}
   \DeclareGraphicsExtensions{,.eps}
\fi

\usepackage[noadjust]{cite}

\begin{document}

\title{Comprehensive Modeling of Three-Phase Distribution Systems via the Bus Admittance Matrix}

\author{Mohammadhafez~Bazrafshan,~\IEEEmembership{Student Member,~IEEE,}
        and~Nikolaos~Gatsis,~\IEEEmembership{Member,~IEEE}% }
 \thanks{The authors are with the Dept. of Electrical \& Computer Engineering,  Univ. of Texas at San Antonio. Emails: \{mohammadhafez.bazrafshan, nikolaos.gatsis\}@utsa.edu. This material is based upon work supported by the National Science Foundation under Grant No. CCF-1421583.}}

\markboth{}%
{}

\maketitle

\begin{abstract}
The theme of this paper is  three-phase distribution system modeling suitable for the Z-Bus load-flow. Detailed models of wye and delta constant-power, constant-current, and constant-impedance loads are presented. Models of transmission lines, step-voltage regulators, and transformers that build the bus admittance matrix (Y-Bus) are  laid out.  The Z-Bus load-flow is then reviewed  and the singularity of the Y-Bus in case of certain transformer connections is rigorously discussed. Based on realistic assumptions and conventional modifications, the invertibility of the Y-Bus  is proved.  Last but not least,  MATLAB scripts that model the components of the IEEE 37-bus, the IEEE 123-bus, the 8500-node feeders,  and the European 906-bus low-voltage feeder are provided. 
\end{abstract}


\section{Introduction}
\IEEEPARstart{A}{n} accurately constructed  bus admittance matrix (Y-Bus) that captures the unbalanced characteristics of distribution networks is the basis of several applications such as \begin{enumerate*}\item three-phase load-flow based on Newton-Raphson \cite{BirtGraffyMcDonalElAbiad1976}, the current injection method \cite{GarciaPereiraCarneiroCostaMartins2000}, or the Z-Bus method \cite{Chen1991pf};
\item three-phase optimal power flow (OPF) using interior-point methods \cite{AraujoPenidoCarneiroPereira2013},  semidefinite relaxations  \cite{DallaneseZhuGiannakis2013}, or successive convex approximations \cite{ZamzamSidiropoulosDallanese2016};
\item voltage security assessment \cite{WuKumagai1982,Overbye1994} through conditions for solution existence \cite{bolognani2016,YuTuritsyn2015,WangBernsteinBoudecPaolone2016,
WangBernsteinBoudecPaolone2016threephase,bazrafshanGatsis2016};
\item optimal system operation by selecting optimal regulator tap settings \cite{RobbinsZhuGarcia2016,BaranFernandes2016} and optimal capacitor switch reconfiguration  \cite{JiangBaldick1996}; and
\item providing real-time voltage solutions by linearizing three-phase power flow equations \cite{BolognaniDorfler2015,KekatosZhangGiannakisBaldick2016,AhmadiMartiMeier2016,Garces2016}.\end{enumerate*}

To facilitate three-phase distribution system studies, such as the above examples,  and to bypass  single-phase simplifications, this paper brings together models---some previously available, some novel---for elements of distribution networks and constructs the Y-Bus matrix.   The Y-Bus matrix incorporates models of three-phase transmission lines, transformers, and step-voltage regulators (SVRs).  

In addition, this paper analyzes how each model affects the invertibility of the Y-Bus matrix and  rigorously proves why a previous proposal in \cite{Gorman1992} for certain transformer connections removes the singularity of the Y-Bus. The invertibility of the Y-Bus is especially important since it allows for computation of load-flow solutions through the Z-Bus method \cite{Chen1991pf} and provides conditions of solution existence  \cite{WangBernsteinBoudecPaolone2016threephase,bazrafshanGatsis2016}.

Three-phase power system modeling is the theme of~\cite{chendillon1974}, although models of voltage-dependent loads and SVRs are not included.  For the forward-backward sweep load-flow (FB-Sweep), ABCD matrices of transmission lines,   distribution transformers, and SVRs are derived in \cite{KerstingBook2001}. Deriving the corresponding admittance matrices from the ABCD parameters is not straightforward---especially for transformers and SVRs.

Three-phase nodal admittance matrices of certain  transformer connections, such as the delta--delta, are singular \cite{chendillon1974, Moorthy2002}.  These singularities are somewhat remedied for the FB-Sweep  in \cite{Kersting1999} by adding independent KCL and KVL equations and in  \cite{Xiao2006} by showing that the zero-sequence components of  voltages do not affect the backward sweep calculations; thus, they render  unique non-zero sequence voltage solutions.  Alternative transformer nodal admittances are suggested in \cite{Dzafic2015} by assembling their symmetrical component circuits.
 
 In the Z-Bus method realm, Ohm's law is  used to  obtain voltage solutions by inverting the Y-Bus at every iteration. To bypass  Y-Bus singularities,  \cite{Gorman1992} advocates adding small admittances from isolated circuits to the ground and \cite{AndersonWollenberg1995} suggests modifying the LDL factorization to use the limit for zero by zero divisions. An alternative approach to avoid Y-Bus singularity is to use equivalent current injections~\cite{Chen1991}.
 
Despite the extensive literature on distribution network modeling \cite{chendillon1974,KerstingBook2001,Moorthy2002,Gorman1992,AndersonWollenberg1995,Chen1991,Xiao2006,Kersting1999,Dzafic2015} and references therein,  precise unified Y-Bus modeling accounting for transmission lines with missing phases and SVRs---as present in IEEE feeders---is lacking. Moreover, there is a deficit in theoretical studies on the invertibility of the Y-Bus and how this is affected by  nodal admittances of various distribution system components.
 
This paper sets out to fill the aforementioned gaps. First, detailed models of wye and delta constant-power, constant-current, and constant-impedance (ZIP) loads,  three-phase transformers, and transmission lines with missing phases  are presented.  A  novel nodal admittance model for SVRs, with and without considering the leakage admittance, is also derived. Specifically, primary-to-secondary  gain matrices for line-to-neutral voltages are derived from first principles. These are combined with a series transmission line, which is the typical configuration in which SVRs are installed in distribution systems. It is worth noting that the  SVR models developed in this work  do not use sequence network matrices (as per \cite{KerstingBook2001}), and can also be utilized in the FB-Sweep algorithm.

Second, this paper precisely explains which distribution system components cause singularities in the Y-Bus and how they can be rectified. In particular, typical SVRs do not render the Y-Bus singular, rather, the Y-Bus singularity is due to self-admittances of certain transformer connections.  Based on linear algebra, this paper shows how connecting small shunt admittances, which has been previously proposed for modification of the Y-Bus~\cite{Gorman1992}, aids in restoring its invertibility.  The Y-Bus invertibility is then proved for  networks that include arbitrary combinations of very general and practical component models---including transmission lines with missing phases---that are found in typical distribution feeders. The invertibility of the Y-Bus is necessary for computing voltage solutions via the Z-Bus method.

Third,  a set of MATLAB scripts is provided online that takes as input the data files for the IEEE 37-bus, the IEEE 123-bus, the 8500-node feeders, and the European 906-bus low voltage feeder \cite{ieeefeederdata}, and models the loads, transmission lines, SVRs, and transformers.   These  scripts further build the corresponding Y-Bus for each feeder and implement the Z-Bus method to compute the voltages.  The power flow solutions obtained are within
$0.75\%$ of benchmark solutions. 

\emph{Paper organization:}  The notation required for three-phase distribution system modeling is introduced in Section \ref{sec:notations}. Models for wye and delta ZIP loads are given in Section \ref{sec:loadmodels}.  Modeling of series elements, i.e., transmission lines, various types of SVRs, and transformers is taken up in Section \ref{sec:serieselements}.  Section \ref{sec:ybuszbus} puts the aforementioned models together, constructs the Y-Bus, reviews the Z-Bus method, and explains the procedures to handle Y-Bus singularities.  A rigorous proof of Y-Bus invertibility is detailed in Section \ref{sec:invertibility}.  Extensions guaranteeing  Y-Bus invertibility under more practical considerations are pursued in Section \ref{sec:extension}.  Section \ref{sec:numtests} provides numerical discussions and load-flow results on distribution test feeders. The paper concludes in Section \ref{sec:conclusion}.

\section{Modeling Notations}
\label{sec:notations}
Power distribution networks comprise two types of elements: 1) shunt elements such as loads and shunt capacitors; and 2) series elements such as three-phase transmission lines, transformers, and SVRs.  Mathematically,  we model a three-phase power distribution network  by an undirected graph $(\mc{N}, \mc{E})$.  The set  $\mc{N}:=\{1,2, \ldots, N\} \cup \{\mr{S}\}$ is the set of nodes and represent the shunt elements, while $\mc{E} := \{(m,n)\} \subseteq \mc{N} \times \mc{N}$ is the set of edges representing the series elements.  Node $\mr{S}$ is considered to be the slack bus connected to the substation. Furthermore, we define the set of neighboring nodes to node $n$ as $\mc{N}_n:= \{ m | (n,m) \in \mc{E} \}$. 

 For a series element, i.e., the edge $(n,m) \in \mc{E}$, let $\Omega_{nm}=\Omega_{mn}$ denote its set of phases.   Define $i_{nm}^{\phi}$ as the current flowing from node $n$ to node $m$ on phase $\phi \in \Omega_{nm}$. Define further the available phases of a node $n$ as $\Omega_n:=\bigcup\nolimits_{m \in \mc{N}_n} \Omega_{nm}$.  Let $\mb{i}_{nm} \in \mbb{C}^{|\Omega_n|}$ collect the currents on all phases flowing from node $n$ to node $m$ such that $\mb{i}_{nm}(\Omega_{nm})=\{i_{nm}^{\phi}\}_{\phi \in \Omega_{nm}}$ and $\mb{i}_{nm}(\{\phi\})=0$ if $\phi \in \Omega_n \setminus \Omega_{nm}$, that is, notation $\mb{i}_{nm}(\Omega_{nm})$  picks the indices of $\mb{i}_{nm}$  that correspond to the phases in $\Omega_{nm}$.
 
 We partition $\mc{N}$ as $\mc{N}=\mc{N}_\mr{Y} \cup  \mc{N}_\Delta \cup \{\mr{S}\} $ where $\mc{N}_{\mr{Y}}$ and $\mc{N}_{\Delta}$ collect wye and delta nodes respectively.   For wye nodes, i.e., $n \in \mc{N}_{\mr{Y}}$, $\Omega_n$ may have one, two, or three available phases.  For delta nodes, i.e., $n \in \mc{N}_{\Delta}$, $\Omega_n$ has at least two available phases, that is $|\Omega_n| \ge 2$.  

For node $n$ and phase $\phi \in \Omega_n$, the complex line to neutral voltage is denoted by $v_n^{\phi}$ and the net current injection is denoted by $i_n^{\phi}$. Moreover, define vectors  $\mb{i}_n=\{i_n^\phi\}_{\phi \in \Omega_n}$, $\mb{v}_n=\{v_n^\phi\}_{\phi \in \Omega_n}$ in $ \mbb{C}^{|\Omega_n|}$ respectively  as the vector of net current injection and complex line to neutral voltages at node $n$. Notice that Ohm's law at each node demands $\mb{i}_n=\sum\nolimits_{m \in \mc{N}_n} \mb{i}_{nm}$. Collect these quantities for all nodes in the vectors $\mb{i}= \{ \mb{i}_n\}_{n \in \mc{N}\backslash \{\mr{S}\}}$  and  $\mb{v}= \{\mb{v}_n\}_{n \in \mc{N} \backslash \{\mr{S} \} }$.

Define further an index set $\mc{J}:=\{1, \ldots, J\}$ where $J=\sum_{n=1}^N |\Omega_n|$, and $j\in \mc{J}$ is a linear index corresponding to a particular pair $(n, \phi)$ with $n \in \mc{N} \backslash \{\mr{S}\}$ and $\phi \in \Omega_n$. In this case,  denote $n=\texttt{Node}[j]$,  and define  $\mc{J}_n:=\{j | \texttt{Node}[j]=n \}$ as the set of linear indices corresponding to node $n$.

\section{Three-Phase ZIP Load Models}
\label{sec:loadmodels}
Due to existence of loads, the nodal net current injection
	$\mb{i}$ is a function of nodal voltages v. This dependence is
	denoted by $\mb{i}_n(\mb{v}_n)$. According to the ZIP load model,   $\mb{i}_n(\mb{v}_n)$  is  composed of  currents from  constant-power loads $\mb{i}_{\mr{PQ}_n} =\bigl\{ i_{\mr{PQ}_n}^{\phi} \bigr\}_{\phi \in \Omega_n}$, constant-current loads $\mb{i}_{\mr{I}_n}= \bigl\{i_{\mr{I}_n}^{\phi}\bigr\}_{\phi \in \Omega_n}$, and constant-impedance loads $\mb{i}_{\mr{Z}_n}=\bigl\{ i_{\mr{Z}_n}^{\phi}\bigr\}_{\phi \in \Omega_n}$.  For  $n \in \mc{N} \backslash \{\mr{S}\}$ and $\phi \in \Omega_n$ we have that
\begin{IEEEeqnarray}{rCl}
i_{n}^\phi(\mb{v}_n) =  i_{\mr{PQ}_n}^{\phi}(\mb{v}_n) + i_{\mr{I}_n}^{\phi} (\mb{v}_n)+ i_{\mr{Z}_n}^{\phi}(\mb{v}_n) \label{eqn:netI}
\end{IEEEeqnarray}
where functions $i_{\mr{PQ}_n}^{\phi}(\mb{v}_n)$, $i_{\mr{I}_n}^{\phi}(\mb{v}_n)$, and $i_{\mr{Z}_n}^\phi(\mb{v}_n)$  are given in Table~\ref{table:wyedeltaloads} for wye and delta connections where  $s_{L_n}^\phi$, $i_{L_n}^{\phi}$, and $y_{L_n}^{\phi}$ are respectively the nominal constant-power,  constant-current, and constant-impedance portions of the ZIP model for wye nodes $n \in \mc{N}_{\mr{Y}}$. Quantities $s_{L_n}^{\phi \phi'}$, $i_{L_n}^{\phi \phi'}$ , and $y_{L_n}^{\phi \phi'}$  are respectively the nominal constant-power,  constant-current, and constant-impedance portions of the ZIP model for nodes $n \in \mc{N}_\Delta$ and over phases $\phi, \phi' \in \Omega_n$. For $n \in \mc{N}_{\Delta}$ and $\phi, \phi' \in \Omega_n$, we have that $s_{L_n}^{\phi \phi'} = s_{L_n}^{\phi' \phi}$, $i_{L_n}^{\phi \phi'}= i_{L_n}^{\phi'\phi}$, and $y_{L_n}^{\phi \phi'}=y_{L_n}^{\phi'\phi}$. Note that a simpler constant-current model is setting $i_{\mr{I}_n}^{\phi}=i_{\mr{L}_n}^{\phi}$. The constant-power or constant-current expressions of Table \ref{table:wyedeltaloads} may also be used to model distributed generation units connected to node $n$  with appropriate signs for $s_{\mr{L}_n}^{\phi}$ and $i_{\mr{L_n}}^{\phi}$.

\begin{table}[t]
\renewcommand{\arraystretch}{1.3} 
\centering
\caption{Wye and delta ZIP load models}
\begin{tabular}{|c|c|c|}
\hline
Current portion & Wye loads & Delta loads \\
\hline
$i_{\mr{PQ}_n}^{\phi}(\mb{v}_n)$ & $-(s_{L_n}^{\phi} / {v_n^\phi})^*$ & $- \sum\limits_{ \phi' \in \Omega_n \backslash \{ \phi \}} \left(\frac{ s_{L_n}^{\phi \phi'}}{v_n^{\phi} - v_n^{\phi'}}\right)^*$ \\
\hline 
$i_{\mr{I}_n}^{\phi}(\mb{v}_n)$  & $ - \frac{v_n^{\phi}}{|v_n^{\phi}|} i_{L_n}^{\phi}$ & $-  \sum\limits_{ \phi' \in \Omega_n \backslash \{ \phi \}} i_{L_n}^{\phi \phi'} \frac{v_n^{\phi} - v_n^{\phi'}}{|v_n^{\phi} - v_n^{\phi'}|}$ \\
\hline
$i_{\mr{Z}_n}^{\phi}(\mb{v}_n)$ & $- y_{L_n}^{\phi}v_n^\phi$ &  $-\sum\limits_{ \phi' \in \Omega_n \backslash \{ \phi \}}   y_{L_n}^{\phi \phi'} (v_n^{\phi} - v_n^{\phi'})$\\
\hline
\end{tabular}
\label{table:wyedeltaloads}
\end{table}

Due to the linear relationship between voltage and currents of constant-impedance loads in Table~\ref{table:wyedeltaloads}, for all nodes $n \in \mc{N}_Y \cup \mc{N}_\Delta$ we have 
\eq{rCl}{
		\mb{i}_{\mr{Z}_n} (\mb{v}_n) &=&  -\mb{Y}_{\mr{L}_n} \mb{v}_n \label{eqn:YLn}
		  } 
where $\mb{Y}_{\mr{L}_n} \in \mbb{C}^{|\Omega_n| \times |\Omega_n|}$ is defined as follows:
\begin{subequations}
\label{eqngroup:YLnspecific}
\eq{rCl}{
	\mb{Y}_{\mr{L}_n} (\phi, \phi) &=&y_{L_n}^{\phi}, n \in \mc{N}_Y, \phi \in \Omega_n \label{eqn:YLnwye} \\
		\mb{Y}_{\mr{L}_n} (\phi, \phi) &=&\sum\limits_{\phi' \in \Omega_n \setminus \{\phi\}} y_{L_n}^{\phi \phi'}, n \in \mc{N}_\Delta, \phi \in \Omega_n \label{eqn:YLndeltadiag}\\
		\mb{Y}_{\mr{L}_n} (\phi, \phi') &=&- y_{L_n}^{\phi \phi'}, n \in \mc{N}_\Delta, \phi \in \Omega_n, \phi' \in \Omega_n \setminus \{\phi\}. \label{eqn:YLndeltaoffdiag}	 \IEEEeqnarraynumspace	
	}
\end{subequations}

\begin{remark}
The ZIP load model only approximately represents the dependencies of nodal injection currents on voltages. For example, the shortcomings of the ZIP model in capturing the high sensitivity of  reactive powers on voltages  is known \cite{loadmodels1995}. In such cases,  more general models, such as ones where  the active and reactive powers are  polynomial or exponential functions of the voltages, can be used \cite[Tables 1.1--1.4]{loadmodels1995biblio}. 

The ZIP load model presented in this section  is in line with traditional distribution system analysis textbooks such as \cite[Chapter 9]{KerstingBook2001}.   However, the dependency of nodal injection currents $\mb{i}_{n}$ on voltages $\mb{v}_n$   is technically imposed by  various  individual   loads (e.g., residential electrical appliances) aggregated at node $n$. Even in the case of steady-state power flow studies,  determining an accurate $\mb{i}_n(\mb{v}_n)$ requires extensive surveys. In such setups, the ZIP load model can  be incorporated to obtain more flexible but yet more complicated models for residential loads; see e.g.,  \cite{McKenna_2016}.   When such complicated load models are utilized, the relationships provided in Table~\ref{table:wyedeltaloads} can be updated accordingly. These modifications will not affect the ensuing discussions on Y-Bus modeling. 
\end{remark}

\section{Modeling of Series Elements}
\label{sec:serieselements}
  Each edge $(n,m)$ represents a series element and is modeled by the following two equations:
\begin{subequations}
\label{eqngroup:genericSeriesForm}
\begin{IEEEeqnarray}{rCl}
\mb{i}_{nm}=\mb{Y}_{nm}^{(n)} \mb{v}_n - \mb{Y}_{nm}^{(m)} \mb{v}_{m} \label{eqn:genericSeriesnside}\\
\mb{i}_{mn}=\mb{Y}_{mn}^{(m)} \mb{v}_{m} - \mb{Y}_{mn}^{(n)} \mb{v}_n. \label{eqn:genericSeriesmside}
\end{IEEEeqnarray}
\end{subequations}
Matrices $\mb{Y}_{nm}^{(n)} \in \mbb{C}^{|\Omega_n|\times |\Omega_n|}$, $\mb{Y}_{nm}^{(m)} \in \mbb{C}^{ |\Omega_n| \times |\Omega_m|}$, $\mb{Y}_{mn}^{(m)} \in \mbb{C}^{ |\Omega_m| \times |\Omega_m|}$, and $\mb{Y}_{mn}^{(n)} \in \mbb{C}^{|\Omega_m| \times |\Omega_n|}$  are determined based on the model of the series element while setting to zero those rows and columns that correspond to missing phases. Matrices $\mb{Y}_{nm}^{(n)}$ and $\mb{Y}_{mn}^{(m)}$ (respectively, $\mb{Y}_{nm}^{(m)}$ and  $\mb{Y}_{mn}^{(n)}$) are called self-admittances (mutual admittances).
\subsection{Transmission lines}
\label{sec:trlines}
\begin{figure}[t]
\centering
\subfloat[]{\includegraphics[scale=0.15]{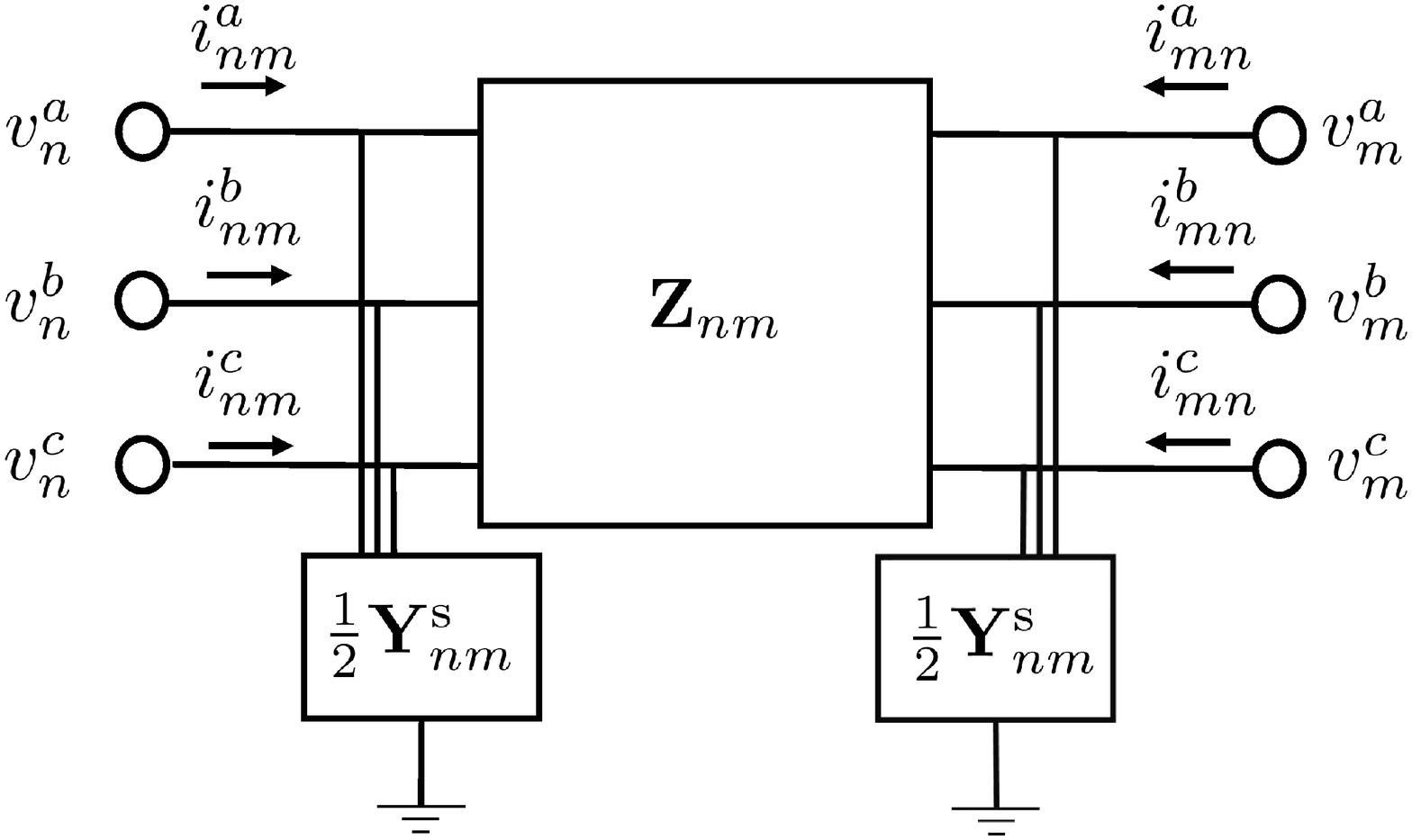} \label{fig:txline}} \quad
\subfloat[ ]{\includegraphics[scale=0.20]{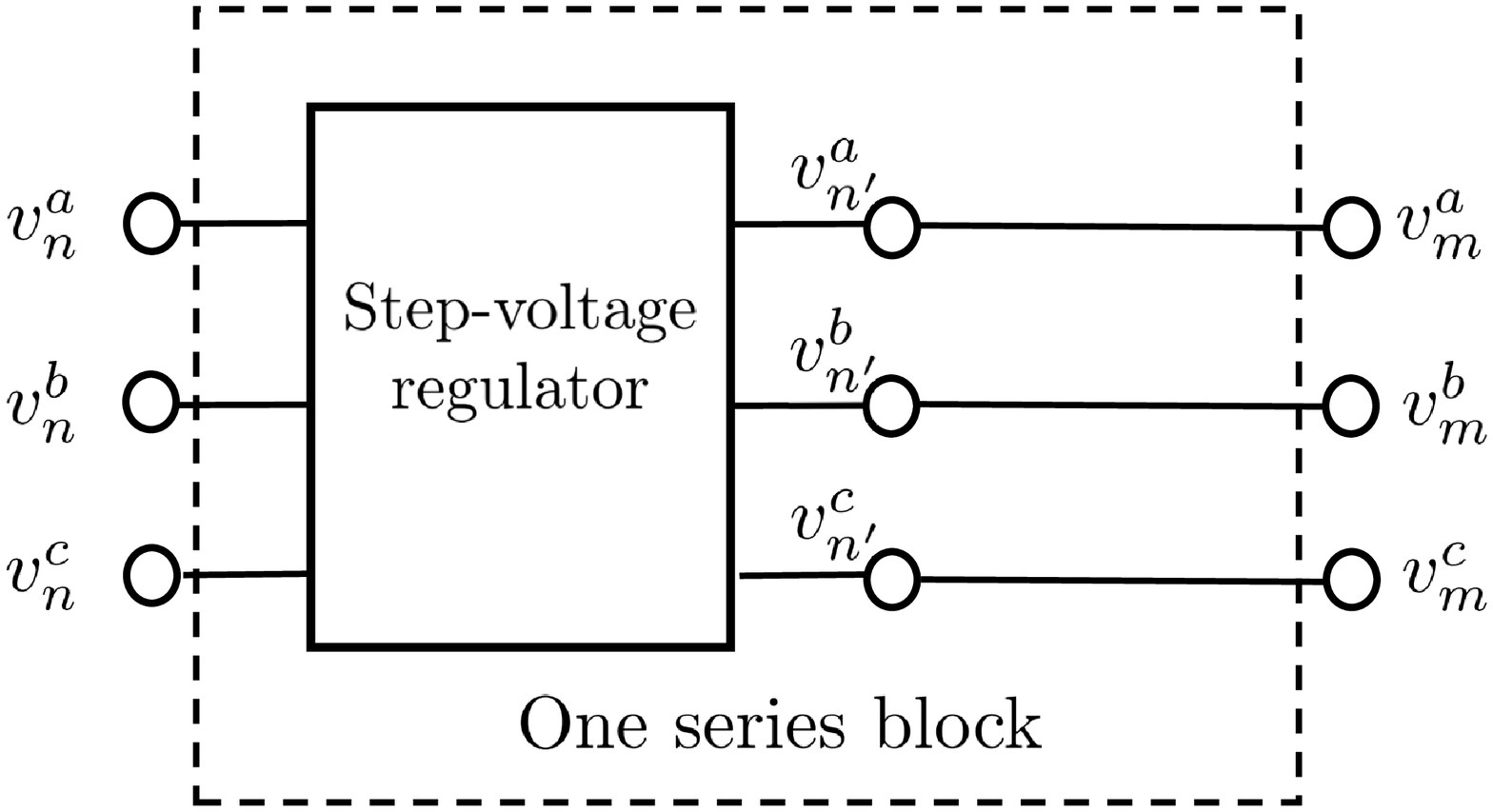}
\label{fig:vr}} 
\caption{\protect\subref{fig:txline} Three-phase transmission line.  \protect\subref{fig:vr} Step-voltage regulator in series with a transmission line.}
\label{fig:trlinevrline}
\end{figure}

 In general, we have that $\Omega_{nm} \subseteq \Omega_n$ for a transmission line on edge $(n,m)$,  and it is possible to have that $|\Omega_{nm}| <|\Omega_{n}|$.  As an example, in the IEEE 13-bus distribution  feeder, the three-phase node 671  is connected to the two-phase node 684 through a two-phase distribution line (config. 604).   The $\pi$-model for transmission lines depicted in Fig.~\ref{fig:trlinevrline}\subref{fig:txline}  yields \cite{KerstingBook2001}
\begin{IEEEeqnarray}{l}
\mb{i}_{nm}(\Omega_{nm})= \left[\frac{1}{2} \mb{Y}_{nm}^{\mr{s}}+ \mb{Z}_{nm}^{-1}\right] \mb{v}_{n}(\Omega_{nm})- \mb{Z}_{nm}^{-1} \mb{v}_m(\Omega_{nm})  \IEEEeqnarraynumspace  \label{eqngroup:txlineactual}
\end{IEEEeqnarray}
where notations $\mb{i}_{nm}(\Omega_{nm})$ and $\mb{v}_n(\Omega_{nm})$ pick the indices of $\mb{i}_{nm}$ and $\mb{v}_n$ that correspond to the phases in $\Omega_{nm}$. For example, if $\mb{v}_n=[v_n^a, v_n^b,v_n^c]^T$ and $\Omega_{nm}=\{b,c\}$, then $\mb{v}_n(\Omega_{nm}) = [v_n^b, v_n^c]^T$.  Matrix $\mb{Z}_{nm}=\mb{Z}_{mn} \in \mbb{C}^{|\Omega_{nm}| \times |\Omega_{nm}|}$ is the series impedance matrix and matrix $\mb{Y}_{nm}^{\mr{s}}= \mb{Y}_{mn}^{\mr{s}} \in \mbb{C}^{|\Omega_{nm}| \times |\Omega_{nm}|}$ is the shunt admittance matrix of the line.  Comparing \eqref{eqngroup:txlineactual} with \eqref{eqngroup:genericSeriesForm} yields
\begin{subequations}
\label{eqngroup:transmissionLineYtildes}
 \begin{IEEEeqnarray}{rClll}
\mb{Y}_{nm}^{(n)}(\Omega_{nm}, \Omega_{nm})&=&\mb{Y}_{mn}^{(m)}(\Omega_{nm}, \Omega_{nm})&=&\frac{1}{2} \mb{Y}_{nm}^{\mr{s}}+ \mb{Z}_{nm}^{-1} \IEEEeqnarraynumspace \\
\mb{Y}_{nm}^{(m)}(\Omega_{nm}, \Omega_{nm})&=&\mb{Y}_{mn}^{(n)}(\Omega_{nm}, \Omega_{nm})&=&\mb{Z}_{nm}^{-1} \IEEEeqnarraynumspace 
 \end{IEEEeqnarray}
 \end{subequations} 
where the notation $\mb{Y}(\Omega_{nm}, \Omega_{nm})$ selects rows and columns of $\mb{Y}$ corresponding to existing phases in $\Omega_{nm}$.

For typical multi-phase transmission lines, the series impedance matrix $\mb{Z}_{nm}$ is symmetric and its real part $\mr{Re}[\mb{Z}_{nm}]$ is positive definite. Based on this property, Lemma \ref{lemma:usefulalgebra} proves that $\mb{Z}_{nm}^{-1}$ exists and that $\mr{Re}[\mb{Z}_{nm}^{-1}]$ is also positive definite.  The shunt admittance matrix $\mb{Y}_{nm}^{\mr{s}}$ is also symmetric and its real part $\mr{Re}[\mb{Y}_{nm}^{\mr{s}}]$ is positive semi-definite.  This specific structure will be shown  to  play a crucial role in guaranteeing the invertibility of a network of three-phase transmission lines. This property will be revisited in Section \ref{sec:invertibility}.

\subsection{Step-voltage regulators}
\label{sec:svrs}
The next series element to be modeled is the SVR, a device which is installed either at the substation or along the feeder to keep nodal voltages within acceptable ranges.  In essence, the SVR is a connection of  auto-transformers  with  adjustable turns ratios.  The turns ratios are dependent on the position of its taps.  The taps are   determined through a control circuit that approximates the voltage drop from the regulator node to the node whose voltage is to be controlled.   In distribution networks, three-phase SVRs are commonly installed in wye, closed-delta, or open-delta configurations \cite{KerstingBook2001}.

The installation of an SVR is in series with a  transmission line similar to Fig.~\ref{fig:trlinevrline}\subref{fig:vr}. The SVR is between nodes $n$ and $n'$ and the transmission line is between nodes $n'$ and $m$. Other than edges $(n',n)$ and $(n',m)$,  no other edges are connected to node $n'$.   With the regulator taps already determined, an equivalent model for the SVR as a single block between nodes $n$ and $m$ can be derived. This equivalent model eliminates node $n'$ and creates a virtual edge $(n,m)$ in place of $(n,n')$.

Voltages and currents at the primary and the secondary of a SVR are  generally related  via
\begin{IEEEeqnarray}{rClrCl}
\mb{v}_n &=& \mb{A}_\mr{v} \mb{v}_{n'}+\mb{Z}_{\mr{R}} \mb{i}_{nn'}, \quad   & 
\mb{i}_{nn'} &=& \mb{A}_\mr{i} \mb{i}_{n'm},  \label{eqngroup:genericReg}
\end{IEEEeqnarray}
where matrices $\mb{A}_\mr{v}, \mb{A}_\mr{i}, \mb{Z}_{\mr{R}} \in \mbb{C}^{\Omega_{n} \times \Omega_{n}}$  will be referred to as voltage gain, current gain, and impedance matrix of the SVR respectively.  The specific entries of the aforementioned matrices $\mb{A}_\mr{v}$,  $\mb{A}_\mr{i}$, and $\mb{Z}_{\mr{R}}$ are determined by the configuration type (wye, closed-delta, or open-delta) as well as the selected tap positions.  A realistic  assumption is that $\Omega_{n} = \Omega_{n'}=\Omega_{nn'}$, which means the primary and secondary of the SVR have the same available phases.   This assumption holds for realistic feeders included in the numerical tests of Section \ref{sec:numtests}.  

It will be shown that, for the three common types of SVRs,  matrix $\mb{Z}_{\mr{R}}$ is diagonal, matrix $\mb{A}_{\mr{v}}$ of \eqref{eqngroup:genericReg}  is invertible, and the following property holds:
\eq{rCl}{\mb{A}_{\mr{v}}^{-1}=\mb{A}_{\mr{i}}^T. \label{eqn:AvAiT}}

Property \eqref{eqn:AvAiT} allows one to obtain nodal admittance models of SVRs that ultimately guarantee the invertibility of the Y-Bus in Section \ref{sec:invertibility}.   Nodal admittance models of SVRs are derived next. From \eqref{eqngroup:genericReg}, by understanding that $\mb{A}_{\mr{v}}^{-1}=\mb{A}_{\mr{i}}^T$, we obtain:
\begin{IEEEeqnarray}{rCl}
\mb{v}_{n'} = \mb{A}_{\mr{i}}^T\left(\mb{v}_{n}-\mb{Z}_{\mr{R}} \mb{A}_{\mr{i}} \mb{i}_{n'm} \right). \label{eqn:genericRegVoltagePrime}
\end{IEEEeqnarray} 
Moreover, the model of \eqref{eqngroup:genericSeriesForm} for the transmission line $(n',m)$ gives the following for the currents $\mb{i}_{n'm}$ and $\mb{i}_{mn'}$: 
\begin{subequations}
\label{eqngroup:transmissionline}
\begin{IEEEeqnarray}{rCl}
\mb{i}_{n'm} &=&\mb{Y}_{n'm}^{(n')}\mb{v}_{n'} - \mb{Y}_{n'm}^{(m)} \mb{v}_m \label{eqn:inprimem} \\
\mb{i}_{mn'} &=& \mb{Y}_{mn'}^{(m)}\mb{v}_m - \mb{Y}_{mn'}^{(n')} \mb{v}_{n'}. \label{eqn:imnprime}
\end{IEEEeqnarray}
\end{subequations}

Using \eqref{eqn:genericRegVoltagePrime} in \eqref{eqn:inprimem} yields 
\eq{rCl}{
\mb{i}_{n'm} = \mb{Y}_{n'm}^{(n')} \mb{A}_{\mr{i}}^T \left( \mb{v}_n - \mb{Z}_{\mr{R}} \mb{A}_{\mr{i}}\mb{i}_{n'm}\right)- \mb{Y}_{n'm}^{(m)} \mb{v}_m.  \label{eqn:in'mreg1}
}
By reorganizing \eqref{eqn:in'mreg1}, we obtain
\eq{rCl}{
 \mb{F}_{\mr{R}} \mb{i}_{n'm} = \mb{Y}_{n'm}^{(n')} \mb{A}_{\mr{i}}^T \mb{v}_n - \mb{Y}_{n'm}^{(m)} \mb{v}_m  \label{eqn:in'mreg} \IEEEeqnarraynumspace }
where   
\eq{rCl}{
\mb{F}_{\mr{R}} &=& \mb{I}_{|\Omega_n|} + \mb{Y}_{n'm}^{(n')}\mb{A}_{\mr{i}}^T \mb{Z}_{\mr{R}} \mb{A}_{\mr{i}} \label{eqn:Freg}
}
and $\mb{I}_{|\Omega_n|}$ is the identity matrix in $\mbb{C}^{|\Omega_n| \times |\Omega_n|}$.   In Lemma \ref{lemma:Fr} we show that  $\mb{F}_{\mr{R}}$ is invertible so that $\mb{i}_{n'm}$ is computed as
\eq{rCl}{
\mb{i}_{n'm} &=& \mb{F}_{\mr{R}}^{-1} \mb{Y}_{n'm}^{(n')} \mb{A}_{\mr{i}}^T \mb{v}_n - \mb{F}_{\mr{R}}^{-1} \mb{Y}_{n'm}^{(m)} \mb{v}_m. \label{eqn:in'mregF}}

Using the equation for the current in \eqref{eqngroup:genericReg}, together with \eqref{eqn:in'mregF}, we write $\mb{i}_{nm}:=\mb{i}_{nn'}$ as a function of $\mb{v}_n$ and $\mb{v}_m$:
\begin{IEEEeqnarray}{rCl}
\mb{i}_{nm} &:=& \mb{i}_{nn'} =  \mb{A}_\mr{i} \mb{i}_{n'm}  \notag \\
&=& \mb{A}_{\mr{i}} \mb{F}_{\mr{R}}^{-1} \mb{Y}_{n'm}^{(n')} \mb{A}_{\mr{i}}^T \mb{v}_n - \mb{A}_{\mr{i}}\mb{F}_{\mr{R}}^{-1} \mb{Y}_{n'm}^{(m)} \mb{v}_m. \IEEEeqnarraynumspace \label{eqn:inmreg}
\end{IEEEeqnarray}
To write $\mb{i}_{mn}$ in terms of $\mb{v}_m$ and $\mb{v}_n$, we first use  \eqref{eqn:genericRegVoltagePrime} in \eqref{eqn:imnprime} to obtain 
\begin{IEEEeqnarray}{rCl}
\mb{i}_{mn}&:=&\mb{i}_{mn'}= \mb{Y}_{mn'}^{(m)}\mb{v}_m -\mb{Y}_{mn'}^{(n')} \mb{v}_{n'} \notag \\
&=& \mb{Y}_{mn'}^{(m)}\mb{v}_m - \mb{Y}_{mn'}^{(n')} \mb{A}_{\mr{i}}^T\left( \mb{v}_n- \mb{Z}_{\mr{R}} \mb{A}_{\mr{i}} \mb{i}_{n'm} \right) \notag \\
&=& \mb{Y}_{mn'}^{(m)} \mb{v}_m - \mb{Y}_{mn'}^{(n')} \mb{A}_{\mr{i}}^T \mb{v}_n \notag \\
&& + \:  \mb{Y}_{mn'}^{(n')} \mb{A}_{\mr{i}}^T\mb{Z}_{\mr{R}} \mb{A}_{\mr{i}} \mb{i}_{n'm}.  \label{eqn:imnreg1} \IEEEeqnarraynumspace
\end{IEEEeqnarray}

Next, we replace $\mb{i}_{n'm}$ in  \eqref{eqn:imnreg1} by its equivalent in \eqref{eqn:in'mregF}:
\eq{rCl}{ 
 \mb{i}_{mn}&=&  \left( \mb{Y}_{mn'}^{(m)} - \mb{Y}_{mn'}^{(n')} \mb{A}_{\mr{i}}^T \mb{Z}_{\mr{R}} \mb{A}_{\mr{i}} \mb{F}_{\mr{R}}^{-1} \mb{Y}_{n'm}^{(m)} \right) \mb{v}_m \notag \\
 \IEEEeqnarraymulticol{3}{r}{-\:  \mb{Y}_{mn'}^{(n')} \left( \mb{I}_{|\Omega_n|} -\mb{A}_{\mr{i}}^T \mb{Z}_{\mr{R}} \mb{A}_{\mr{i}} \mb{F}_{\mr{R}}^{-1} \mb{Y}_{n'm}^{(n')}\right) \mb{A}_{\mr{i}}^T  \mb{v}_n. \label{eqn:imnreg} \IEEEeqnarraynumspace}
}

Equations \eqref{eqn:inmreg} and \eqref{eqn:imnreg} conform to the model of~\eqref{eqngroup:genericSeriesForm}  and yield the following matrices for the SVRs:
\begin{subequations}
	\label{eqngroup:SVRnonIdealYtildes}
	\eq{rCl}{
	\mb{Y}_{nm}^{(n)} &=&  \mb{A}_{\mr{i}} \mb{F}_{\mr{R}}^{-1} \mb{Y}_{n'm}^{(n')}  \mb{A}_{\mr{i}}^T   \label{eqn:SVRnonIdealYnmn} \\
	\mb{Y}_{nm}^{(m)} &=&  \mb{A}_{\mr{i}}\mb{F}_{\mr{R}}^{-1} \mb{Y}_{n'm}^{(m)} \label{eqn:SVRnonIdealYnmm} \\
	\mb{Y}_{mn}^{(m)} &=& \mb{Y}_{mn'}^{(m)} - \mb{Y}_{mn'}^{(n')}  \mb{A}_{\mr{i}}^T  \mb{Z}_{\mr{R}} \mb{A}_{\mr{i}} \mb{F}_{\mr{R}}^{-1} \mb{Y}_{n'm}^{(m)} \label{eqn:SVRnonIdealYmnm} \\
	\mb{Y}_{mn}^{(n)} &=& \mb{Y}_{mn'}^{(n')}  \mb{F}_{\mr{R}}^{-T} \mb{A}_{\mr{i}}^T  \label{eqn:SVRnonIdealYmnn} \IEEEeqnarraynumspace
}
\end{subequations}
where in \eqref{eqn:SVRnonIdealYmnn} we have used the identity $\mb{F}_{\mr{R}}^{-T}= \mb{I}_{|\Omega_n|} -\mb{A}_{\mr{i}}^T \mb{Z}_{\mr{R}} \mb{A}_{\mr{i}} \mb{F}_{\mr{R}}^{-1} \mb{Y}_{n'm}^{(n')}$ provided from Lemma \ref{lemma:Fr}.

In most cases, it turns  out that the per unit series impedance of SVRs significantly depends on  the tap position and is zero, rendering $\mb{Z}_{\mr{R}}=\mb{O}$ and $\mb{F}_{\mr{R}}=\mb{I}_{|\Omega_n|}$. Thus, we can obtain the following matrices for SVRs with ideal auto-transformers:
\begin{subequations}
\label{eqngroup:SVRIdealYtildes}
\begin{IEEEeqnarray}{rclrcl}
\mb{Y}_{nm}^{(n)} &=& \mb{A}_{\mr{i}}  \mb{Y}_{n'm}^{(n')} \mb{A}_{\mr{i}}^T, \quad 
&\mb{Y}_{nm}^{(m)} &=& \mb{A}_{\mr{i}}\mb{Y}_{n'm}^{(m)}, \label{eqn:vrynmm} \\
\mb{Y}_{mn}^{(m)} &=&\mb{Y}_{mn'}^{(m)}, \quad
&\mb{Y}_{mn}^{(n)} &=& \mb{Y}_{mn'}^{(n')} \mb{A}_\mr{i}^T.  \label{eqn:vrymnn}
\end{IEEEeqnarray}
\end{subequations} 
Notice in \eqref{eqngroup:SVRnonIdealYtildes} and \eqref{eqngroup:SVRIdealYtildes} that the matrix sizes conform due to the initial assumption that $\Omega_{n}=\Omega_{n'}$.

\begin{figure}[t]
\centering
\includegraphics[scale=0.22]{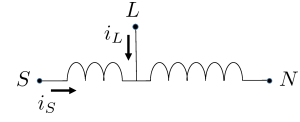}
\caption{Model of a single-phase auto-transformer. For brevity, the series impedance is not drawn.}
\label{fig:type-b-single-regulator}
\vspace{-0.3cm}
\end{figure}

To derive the gain matrices  $\mb{A}_{\mr{v}}$ and $\mb{A}_{\mr{i}}$ as well as the impedance matrix $\mb{Z}_{\mr{R}}$ for the three SVR configurations (wye, closed-delta, and open-delta), it is essential to grasp the basic model of single-phase auto-transformers. A diagram of a single-phase auto-transformer  is given in Fig.~\ref{fig:type-b-single-regulator}.\footnote{Depending on the connection of the shunt winding of auto-transformers,  SVRs can either be of type A or type B. Since type B SVRs are more common in distribution networks \cite{KerstingBook2001}, the regulator models developed correspond only to this type of SVRs.   The same procedure can be followed through to obtain models for type A SVRs. This has been done for the 8500-node feeder of Section \ref{sec:numtests}. }  Given the taps, the effective regulator ratio $a_R$ is determined as follows:
\begin{IEEEeqnarray}{rCl}
a_R = 1 \mp 0.00625 \texttt{tap}. \label{eqn:1phaseregtap}
\end{IEEEeqnarray}
 Having determined $a_R$,  the voltages and currents on the two sides of the auto-transformer relate as follows:
 \begin{subequations}
\label{eqngroup:1phasereg}
\begin{IEEEeqnarray}{rCl}
v_{S} - v_N &=& a_R (v_L - v_N)+z_R i_S \label{eqn:1phaseregv} \\
 i_S &=& -\frac{1}{a_R} i_L  \label{eqn:1phaseregi}
\end{IEEEeqnarray}
\end{subequations}
where $z_R$ is the series impedance of the auto-transformer. 
The voltage and current relationships in \eqref{eqngroup:1phasereg} are leveraged  to derive the gain matrices $\mb{A}_{\mr{v}}$ and $\mb{A}_{\mr{i}}$ for three  types of SVR connections. Table~\ref{table:vrmatrices} summarizes the resulting formulas.  Notice $\mb{A}_{\mr{v}}$ is invertible and that $\mb{A}_{\mr{v}}^{-1}= \mb{A}_{\mr{i}}^T$ holds for all SVRs; a crucial property for guaranteeing Y-Bus invertibility.  
\begin{table*}[t]
\centering
\renewcommand{\arraystretch}{1.75}
\caption{Voltage gain, current gain, and impedance matrices for the three common configurations  of step-voltage regulators}
\begin{tabular}{|c|c|c|c|} 
\hline 
SVR connection & wye-connected &  closed-delta  & open-delta \\
\hline 
Voltage gain $\mb{A}_{\mr{v}}$ & $\bmat{ a_{R_a} & 0 & 0 \\ 0 & a_{R_b} & 0 \\ 0 & 0 & a_{R_c}}$ &  $\bmat{ a_{R_{ab}}  & 1- a_{R_{ab}} & 0 \\  0 & a_{R_{bc}} & 1- a_{R_{bc}} \\  1- a_{R_{ca}} & 0  & a_{R_{ca}}}$  & $\bmat{ a_{R_{ab}} & 1- a_{R_{ab}} & 0 \\ 0 & 1 &  0 \\ 0 & 1-a_{R_{cb}} & a_{R_{cb}}}$ \\
\hline 
Current gain $\mb{A}_{\mr{i}}$ & $\bmat{\frac{1}{a_{R_a}} & 0 & 0 \\ 0 & \frac{1}{a_{R_b}} & 0 \\ 0 & 0 & \frac{1}{a_{R_c}}}$ & $ \bmat{  a_ {R_{ab}} & 0 & 1-a_{R_{ca}} \\ 1 - a_{R_{ab}} & a_{R_{bc}} & 0 \\  0 & 1- a_{R_{bc}} & a_{R_{ca}} } ^{-1}$ &  $\bmat{ \frac{1}{a_{R_{ab}}} & 0 & 0 \\ 1-\frac{1}{a_{R_{ab}}} & 1 & 1-\frac{1}{a_{R_{cb}}} \\ 0 & 0 & \frac{1}{a_{R_{cb}}}}$ \\
\hline 
Impedance matrix $\mb{Z}_{\mr{R}}$ &  $\bmat{ z_{R_{a}} & 0 & 0 \\ 0 & z_{R_b} & 0 \\ 0 & 0 & z_{R_c}}$ & $\bmat{ z_{R_{ab}} & 0 & 0 \\ 0 & z_{R_{bc}} & 0 \\ 0 & 0 & z_{R_{ca}}}$ &  $\bmat{ z_{R_{ab}} & 0 & 0 \\ 0 & 0 & 0 \\ 0 & 0 & z_{R_{cb}}}$ \\
\hline
\end{tabular}
\label{table:vrmatrices}
\end{table*}

\begin{figure*}[t]
\centering
\subfloat[]{\includegraphics[scale=0.15]{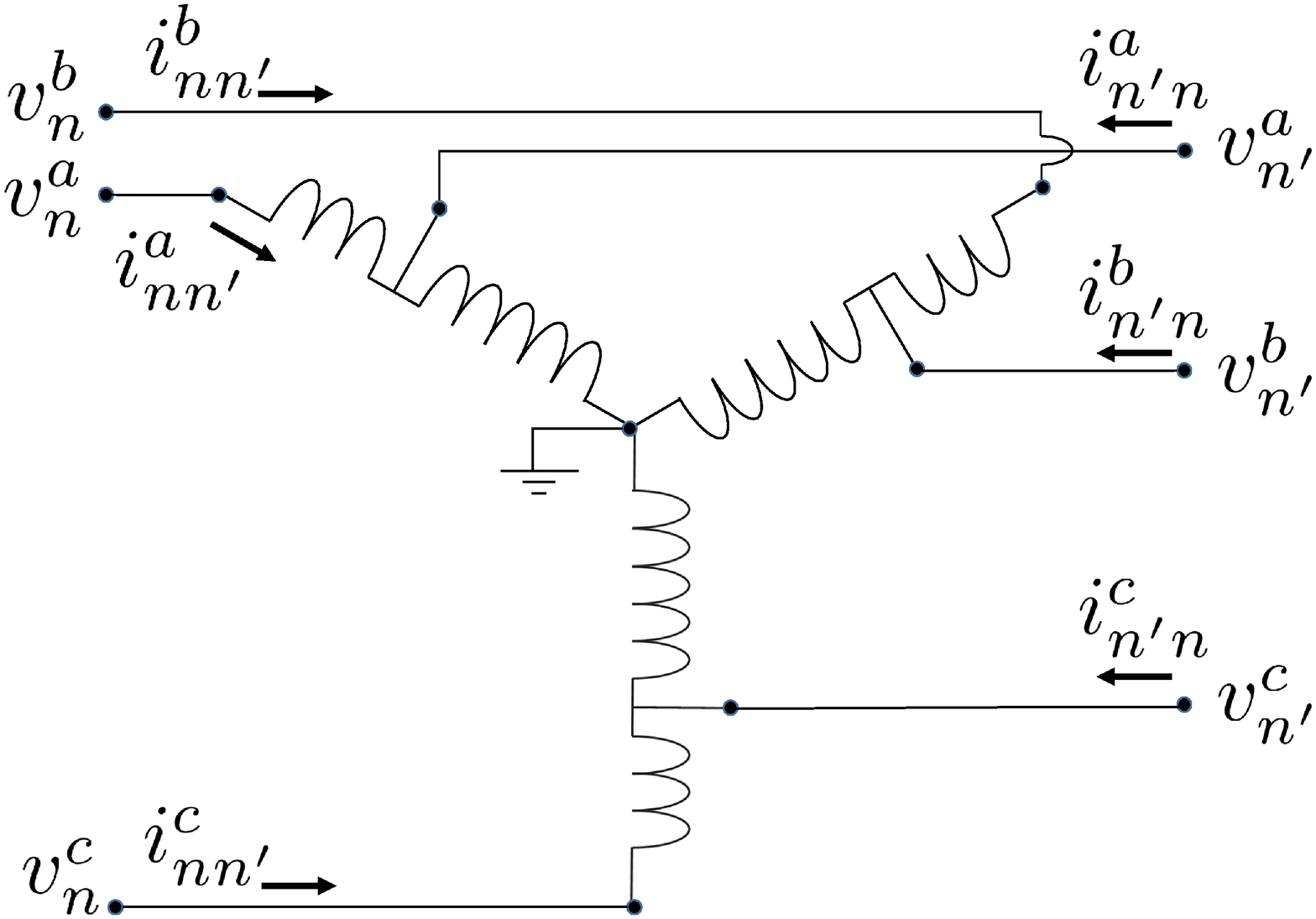} \label{fig:wye-connected regulator}} 
\hspace{1cm}
\subfloat[]{\includegraphics[scale=0.15]{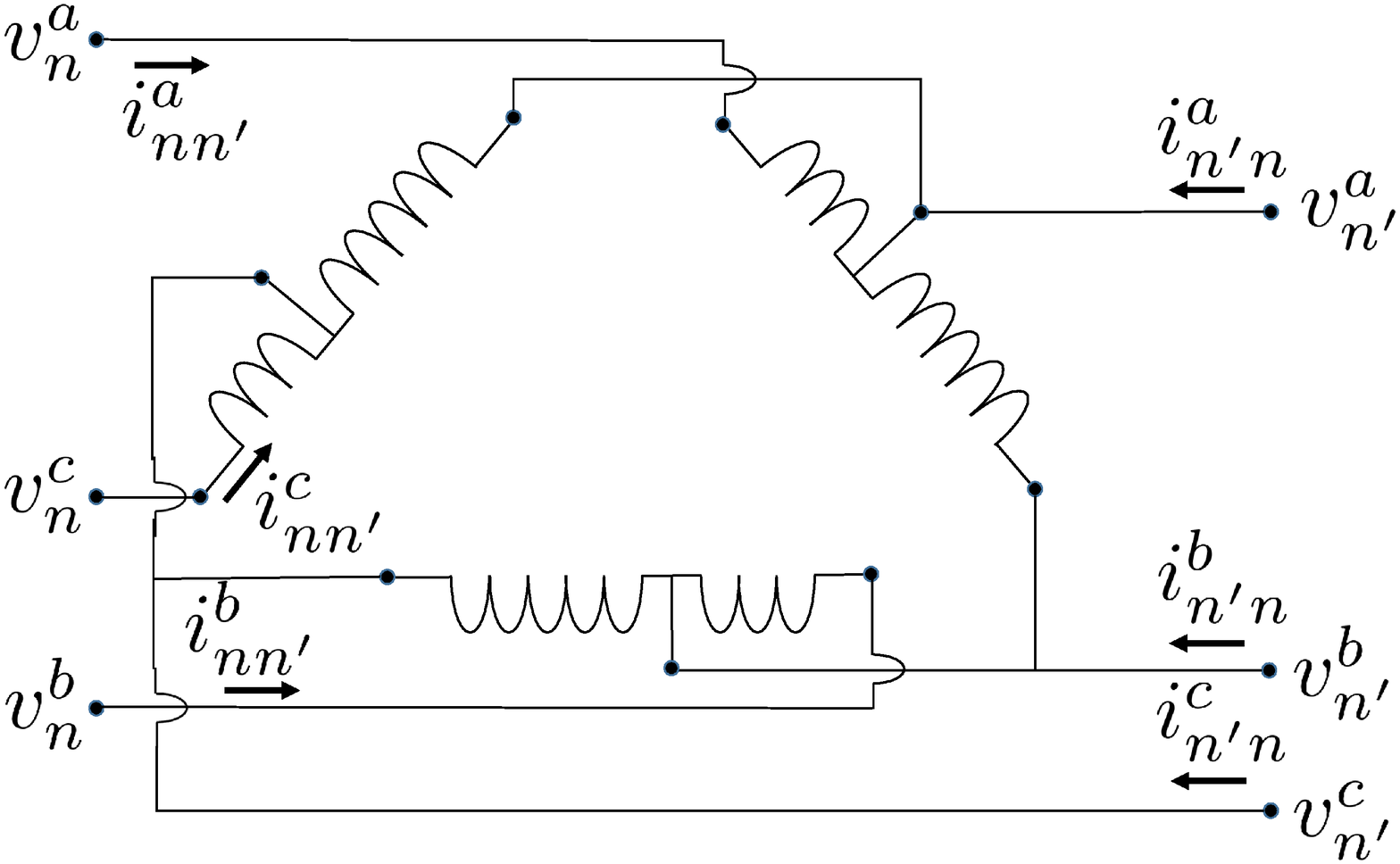}\label{fig:closed-delta-regulator}} \hspace{1cm}
\subfloat[]{\includegraphics[scale=0.15]{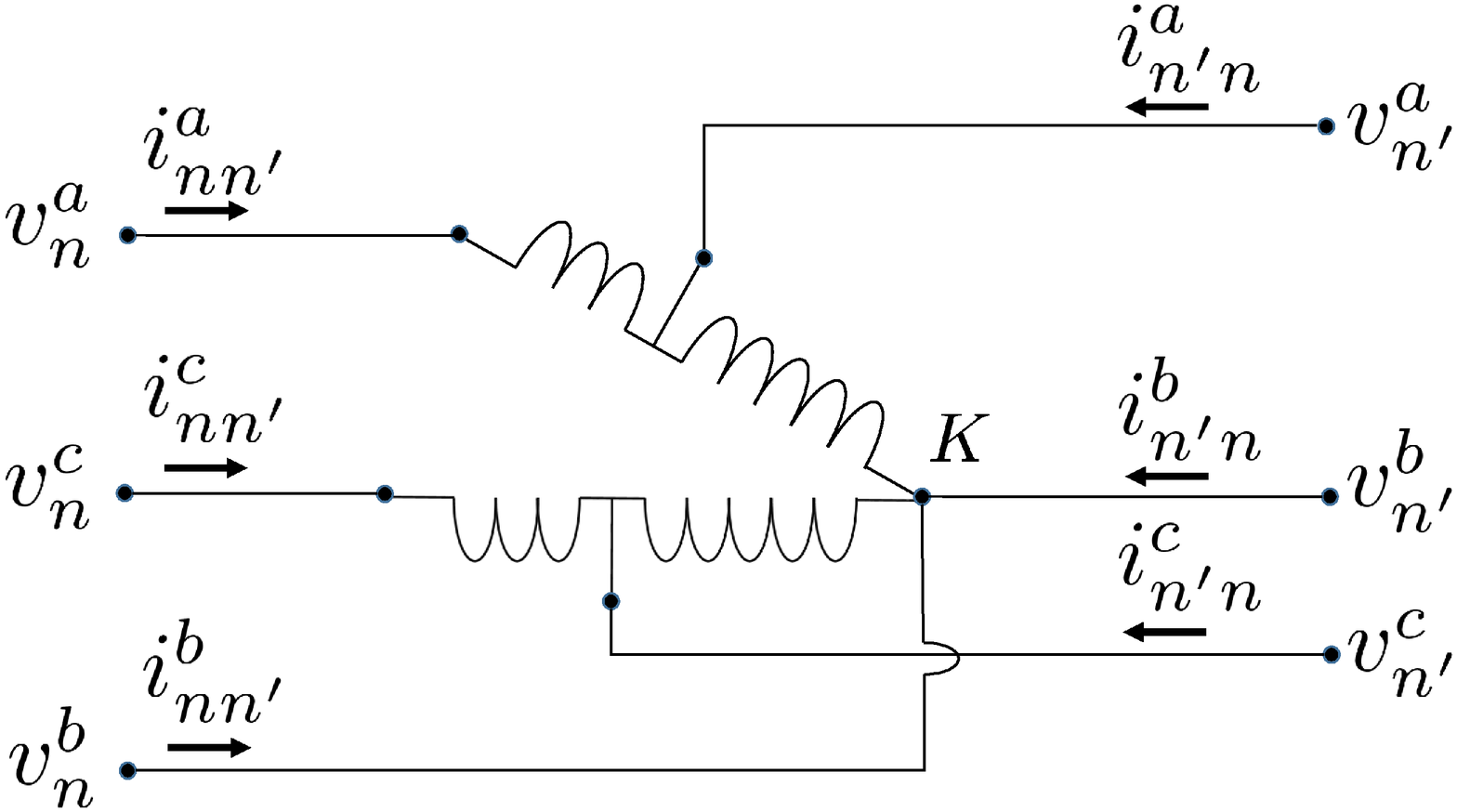} \label{fig:open-delta-regulator}}
\caption{ Three types of step-voltage regulators; \protect\subref{fig:wye-connected regulator} wye-connected \protect\subref{fig:closed-delta-regulator} closed-delta\protect\subref{fig:open-delta-regulator} open-delta. For brevity, the series impedances are not drawn.}
\label{fig:vrmodels}
\end{figure*}

\subsubsection{Wye-connected SVRs}
The diagram of a wye-connected SVR is given in Fig.~\ref{fig:vrmodels}\subref{fig:wye-connected regulator}.   Using  \eqref{eqn:1phaseregv} for  the three auto-transformers in Fig.~\ref{fig:vrmodels}\subref{fig:wye-connected regulator}  yields the following equations:
\begin{subequations}
	\label{eqngroup:wye-svr-voltages}
\begin{IEEEeqnarray}{rCl}
v_n^a &=& a_{R_a} v_{n'}^a + z_{R_a} i_{nn'}^a \label{eqn:wye-svrva}  \\
  v_n^b &=& a_{R_b} v_{n'}^b + z_{R_b} i_{nn'}^b \label{eqn:wye-svrvb} \\
   v_n^c &=& a_{R_c} v_{n'}^c+ z_{R_c} i_{nn'}^c  \label{eqn:wye-svrc} \IEEEeqnarraynumspace
   \end{IEEEeqnarray}
\end{subequations}
where $a_{R_a}$, $a_{R_b}$, and $a_{R_c}$ are the effective regulator ratios  determined by the corresponding taps at phases $a$, $b$, and $c$ following \eqref{eqn:1phaseregtap}. Similarly, $z_{R_a}$, $z_{R_b}$, and $z_{R_c}$ are per phase series impedances of each SVR.   From \eqref{eqngroup:wye-svr-voltages}, one obtains the SVR voltage gain $\mb{A}_{\mr{v}}$ and impedance matrix $\mb{Z}_{\mr{R}}$: 
\eq{rCl}{
\mb{A}_\mr{v}= \bmat{ a_{R_a} & 0 & 0 \\ 0 & a_{R_b} & 0 \\ 0 & 0 & a_{R_c}}, \mb{Z}_{\mr{R}} = \bmat{ z_{R_a} & 0 & 0 \\ 0 & z_{R_b} & 0 \\ 0 & 0 & z_{R_c}}. \label{eqngroup:AvZ-wye-svr} \IEEEeqnarraynumspace
}
From \eqref{eqn:1phaseregi}, we obtain 
\eq{rCl}{
i_{nn'}^a = -\frac{1}{a_{R_a}} i_{n'n}^a, 
i_{nn'}^b = -\frac{1}{a_{R_b}} i_{n'n}^b, 
i_{nn'}^c = -\frac{1}{a_{R_c}} i_{n'n}^c.  \IEEEeqnarraynumspace \label{eqngroup:wye-svr-currents}}
Using \eqref{eqngroup:wye-svr-currents} and recognizing that $\mb{i}_{n'n} = - \mb{i}_{n'm}$ holds in \eqref{eqngroup:genericReg}, we compute $\mb{A}_{\mr{i}}$
\begin{IEEEeqnarray}{rCl}
\mb{A}_{\mr{i}} = \bmat{\frac{1}{a_{R_a}} & 0 & 0 \\ 0 & \frac{1}{a_{R_b}} & 0 \\ 0 & 0 & \frac{1}{a_{R_c}}}. \IEEEeqnarraynumspace \label{eqn:Ai-wye-svr}
\end{IEEEeqnarray}

\subsubsection{Closed-delta connected SVRs}
The diagram of closed-delta connected SVR is  drawn in Fig.~\ref{fig:vrmodels}\subref{fig:closed-delta-regulator}.   Leveraging \eqref{eqngroup:1phasereg} for the three auto-transformers between phases $ab$, $bc$ and $ca$ at node $n'$, the voltage equations are
\begin{subequations}
\label{eqngroup:closed-delta-svr-voltages}
\begin{IEEEeqnarray}{rCl}
v_n^a - v_{n'}^b &=& a_{R_{ab}} (v_{n'}^a-v_{n'}^b) + z_{R_{ab}} i_{nn'}^a \notag \\ 
\Rightarrow v_n^a &=& a_{R_{ab}} v_{n'}^a + (1- a_{R_{ab}}) v_{n'}^b +z_{R_{ab}} i_{nn'}^a\label{eqn:closed-delta-svr-va}  \IEEEeqnarraynumspace\\
v_{n}^b- v_{n'}^c&=& a_{R_{bc}} (v_{n'}^b-v_{n'}^c)+z_{R_{bc}}i_{nn'}^b \notag \\
\Rightarrow v_n^b &=& a_{R_{bc}} v_{n'}^b + (1- a_{R_{bc}}) v_{n'}^c  +z_{R_{bc}}i_{nn'}^b\label{eqn:closed-delta-svr-vb} \IEEEeqnarraynumspace \\
v_{n}^c - v_{n'}^a &=& a_{R_{ca}} (v_{n'}^c- v_{n'}^a)+z_{R_{ca}}i_{nn'}^c \notag \\
\Rightarrow v_{n}^c &=& a_{R_{ca}} v_{n'}^c + (1- a_{R_{ca}}) v_{n'}^a+z_{R_{ca}}i_{nn'}^c \label{eqn:closed-delta-svr-vc} \IEEEeqnarraynumspace
\end{IEEEeqnarray}
\end{subequations}
where $a_{R_{ab}}$,$a_{R_{bc}}$, and $a_{R_{ca}}$ are the effective regulator ratios and $z_{R_{ab}}$, $z_{R_{bc}}$, and $z_{R_{ca}}$ are series impedances on phases $ab$, $bc$, and $ca$ at node $n'$ respectively.   From \eqref{eqngroup:closed-delta-svr-voltages}, we find
\begin{subequations}
	\label{eqngroup:AvZ-closed-delta-svr}
\begin{IEEEeqnarray}{rCl}
\mb{A}_{\mr{v}} &=& \bmat{ a_{R_{ab}}  & 1- a_{R_{ab}} & 0 \\  0 & a_{R_{bc}} & 1- a_{R_{bc}} \\  1- a_{R_{ca}} & 0  & a_{R_{ca}}} \label{eqn:Av-closed-delta-svr} \\
\mb{Z}_{\mr{R}} &=& \bmat{ z_{R_{ab}} & 0 & 0 \\ 0 & z_{R_{bc}} & 0 \\ 0 & 0 & z_{R_{ca}}}. \label{eqn:Z-closed-delta-svr}
\end{IEEEeqnarray}
\end{subequations}
Writing KCL for node $n'$ at phases $a,b,c$ and considering the current equation in \eqref{eqn:1phaseregi} for auto-transformers yields:
\begin{subequations}
\label{eqngroup:closed-delta-regi}
\begin{IEEEeqnarray}{rCl}
i_{n'n}^a &=& - a_{R_{ab}} i_{nn'}^a - (1-a_{R_{ca}}) i_{nn'}^c \label{eqn:closed-delta-reginprimena} \\
i_{n'n}^b &=& - a_{R_{bc}} i_{nn'}^b - (1- a_{R_{ab}}) i_{nn'}^a \label{eqn:closed-delta-reginprimenb} \\
i_{n'n}^c &=& -a_{R_{ca}} i_{nn'}^c - (1-a_{R_{bc}}) i_{nn'}^b. \label{eqn:closed-delta-reginprimenc}
\end{IEEEeqnarray}
\end{subequations}
From \eqref{eqngroup:closed-delta-regi} and using $\mb{i}_{n'n} = - \mb{i}_{n'm}$ in \eqref{eqngroup:genericReg}, we find
\begin{IEEEeqnarray}{rCl}
\mb{A}_{\mr{i}} = \bmat{  a_ {R_{ab}} & 0 & 1-a_{R_{ca}} \\ 1 - a_{R_{ab}} & a_{R_{bc}} & 0 \\  0 & 1- a_{R_{bc}} & a_{R_{ca}} } ^{-1}.  \label{eqn:closed-delta-regAi}
\end{IEEEeqnarray}

\subsubsection{Open-delta connected SVRs}
The circuit diagram of the open-delta SVR is given in Fig.~\ref{fig:vrmodels}\subref{fig:open-delta-regulator}. In Fig.~\ref{fig:vrmodels}\subref{fig:open-delta-regulator}, using the voltage equations in \eqref{eqngroup:1phasereg} and noticing that $v_{n}^b= v_{n'}^b$ (phase $b$ of node $n$ is connected to phase $b$ of node $n'$) yields
\begin{subequations}
\label{eqngroup:open-delta-svr-voltages}
\begin{IEEEeqnarray}{rCl}
v_{n}^a - v_{n'}^b &=& a_{R_{ab}} (v_{n'}^a - v_{n'}^b)+z_{R_{ab}} i_{nn'}^a \notag \\
 \Rightarrow v_n^a &=& a_{R_{ab}} v_{n'}^a + (1- a_{R_{ab}}) v_{n'}^b+z_{R_{ab}} i_{nn'}^a \label{eqn:open-delta-svr-va} \\
v_{n}^c - v_{n'}^b &=& a_{R_{cb}}(v_{n'}^c - v_{n'}^b) + z_{R_{cb}} i_{nn'}^c \notag \\
 \Rightarrow v_{n}^c &=& a_{R_{cb}} v_{n'}^c + (1- a_{R_{cb}}) v_{n'}^b+z_{R_{cb}}i_{nn'}^c \label{eqn:open-delta-svr-vc}
\end{IEEEeqnarray}
\end{subequations}
where $a_{R_{ab}}$ and $a_{R_{cb}}$ are the effective regulator ratios  and $z_{R_{ab}}$ and $z_{R_{cb}}$ are the series impedances on phases  $ab$ and $cb$ of node $n$.  Using \eqref{eqngroup:open-delta-svr-voltages} and the equality $v_{n}^b= v_{n'}^b$, the gain matrix $\mb{A}_{\mr{v}}$, and the impedance matrix $\mb{Z}_{\mr{R}}$  are obtained as follows
\begin{subequations}
\begin{IEEEeqnarray}{rCl}
\mb{A}_{\mr{v}} &=& \bmat{ a_{R_{ab}} & 1- a_{R_{ab}} & 0 \\ 0 & 1 &  0 \\ 0 & 1-a_{R_{cb}} & a_{R_{cb}}} \label{eqn:Av-open-delta-svr} \\
\mb{Z}_{\mr{R}}&=& \bmat{ z_{R_{ab}} & 0 & 0 \\ 0 & 0 & 0 \\ 0 & 0 & z_{R_{cb}}}. \label{eqn:Z-open-delta-svr}
\end{IEEEeqnarray}
\end{subequations}
For the two SVRs in Fig.~\ref{fig:vrmodels}\subref{fig:open-delta-regulator}, from \eqref{eqn:1phaseregi} we obtain
\begin{IEEEeqnarray}{rClrCl}
i_{nn'}^a &=& - \frac{1}{a_{R_{ab}}} i_{n'n}^a, \quad &
i_{nn'}^c & =& -\frac{1}{a_{R_{cb}}} i_{n'n}^c. \label{eqn:open-delta-regi}
\end{IEEEeqnarray}
Applying KCL at node $K$ of Fig.~\ref{fig:vrmodels}\subref{fig:open-delta-regulator} and subsequently introducing \eqref{eqn:open-delta-regi} yields $i_{nn'}^{b}$ as follows:
\begin{align}
i_{nn'}^{b} & = - i_{nn'}^{a} - i_{nn'}^{c} - i_{n'n}^{a} - i_{n'n}^{b} - i_{n'n}^{c} \notag\\
& = \left(\frac{1}{a_{R_{ab}}}  -1\right) i_{n'n}^{a} - i_{n'n}^{b} -  \left(\frac{1}{a_{R_{cb}}}  -1\right)  i_{n'n}^{c}.
\label{eqn:open-delta-regib}
\end{align}
Combining $\mb{i}_{n'n}= - \mb{i}_{n'm}$ with~\eqref{eqn:open-delta-regi} and~\eqref{eqn:open-delta-regib} yields $\mb{A}_{\mr{i}}$:
\begin{IEEEeqnarray}{rCl}
\mb{A}_{\mr{i}} = \bmat{ \frac{1}{a_{R_{ab}}} & 0 & 0 \\ 1-\frac{1}{a_{R_{ab}}} & 1 & 1-\frac{1}{a_{R_{cb}}} \\ 0 & 0 & \frac{1}{a_{R_{cb}}}}.
\end{IEEEeqnarray}

\subsection{Three-phase transformers} 
\label{sec:transformers}
In distribution systems, a three-phase transformer can be appropriately represented by two blocks, namely, a series block representing the per unit leakage admittance, and a shunt block modeling transformer core losses. Transformer core losses can technically be treated as voltage-dependent loads similar to the ones in Table~\ref{table:wyedeltaloads},  with functions of  $\mb{v}$  given in \cite{Chen1991}; notice though that this information may  not be readily available for many distribution feeders.  The most common distribution transformers   include delta--wye grounded, wye--delta, wye--wye, open-wye--open-delta, delta--delta, and open-delta--open-delta for which the ABCD parameters are listed in \cite{KerstingBook2001}.  The corresponding nodal admittances of these connections are given in \cite{Chen1991} and \cite{Chen1992} and are listed in Table~\ref{table:transformerConnections} where 

{\small \begin{subequations}
\label{eqngroup:y1y2y3}
\begin{IEEEeqnarray}{rCl}
\mb{Y}_1&=&\bmat{y_t & 0 & 0  \\ 0 & y_t & 0\\ 0 & 0 & y_t},
\mb{Y}_2= \frac{1}{3} \bmat{ 2y_t & -y_t & -y_t \\ -y_t & 2y_t & -y_t \\ -y_t & -y_t & 2y_t}\\
\mb{Y}_3 &=& \frac{1}{\sqrt{3}} \bmat{ -y_t & y_t & 0 \\ 0 & -y_t & y_t \\ y_t & 0 & -y_t},
\mb{Y}_4= \frac{1}{3}\bmat{ y_t & -y_t & 0 \\ -y_t & 2y_t & -y_t \\ 0 & -y_t & y_t} \IEEEeqnarraynumspace \\
\mb{Y}_5&=& \bmat{y_t & 0 \\ 0 & y_t}, \mb{Y}_6= \frac{1}{\sqrt{3}} \bmat{-y_t & y_t & 0 \\ 0 & -y_t & y_t}. \IEEEeqnarraynumspace
\end{IEEEeqnarray}
\end{subequations}}
\noindent and $y_t$ is the per unit leakage admittance. In the next section, the bus admittance matrix Y-Bus is formulated. It turns out that rank-deficiency of matrices other than $\mb{Y}_1$ in \eqref{eqngroup:y1y2y3} are the only source of Y-Bus singularity. A  method is presented to remedy this issue.
\begin{table}
\renewcommand{\arraystretch}{1.3} 
\centering
\caption{Admittance matrices of the most common transformer connections in distribution systems.}
\begin{tabular}{|cc|c|c|c|c|}
\hline
\multicolumn{2}{|c|}{Transformer connection} & \multicolumn{4}{|c|}{Matrices} \\
\hline
Node $n$  & Node $m$  & $\mb{Y}_{nm}^{(n)}$  & $\mb{Y}_{nm}^{(m)}$  & $\mb{Y}_{mn}^{(m)}$ & $\mb{Y}_{mn}^{(n)}$  \\ 
\hline
Wye-G & Wye-G & $\mb{Y}_1$ & $\mb{Y}_1$ & $\mb{Y}_1$ & $\mb{Y}_1$ \\
\hline 
Wye-G & Wye & $\mb{Y}_2$ & $\mb{Y}_2$ & $\mb{Y}_2$ & $\mb{Y}_2$ \\
\hline 
Wye-G & Delta & $\mb{Y}_1$ & $-\mb{Y}_3$ & $\mb{Y}_2$ & $-\mb{Y}_3^T$ \\
\hline 
Wye  & Wye & $\mb{Y}_2$ & $\mb{Y}_2$ & $\mb{Y}_2$ & $\mb{Y}_2$ \\
\hline 
Wye & Delta & $\mb{Y}_2$ & $-\mb{Y}_3$ & $\mb{Y}_2$ & $-\mb{Y}_3^T$ \\
\hline 
Delta & Delta & $\mb{Y}_2$ & $\mb{Y}_2$ & $\mb{Y}_2$ & $\mb{Y}_2$ \\
\hline
Open-Delta & Open-Delta & $\mb{Y}_4$ & $\mb{Y}_4$ & $\mb{Y}_4$ & $\mb{Y}_4$ \\
\hline 
Open-Wye & Open-Delta & $\mb{Y}_5$ & $-\mb{Y}_6$ & $ \mb{Y}_4$ & $-\mb{Y}_6^T$ \\
\hline
\end{tabular}
\label{table:transformerConnections}
\end{table}

\section{Y-Bus Construction and the Z-Bus Method}
\label{sec:ybuszbus}
The application of KCL at each node $n \in \mc{N}$ leads to
\begin{IEEEeqnarray}{rCl}
\mb{i}_n (\mb{v}_n) = \sum\limits_{m \in \mc{N}_n} \mb{i}_{nm},  n \in \mc{N}.  \label{eqn:KCL}
\end{IEEEeqnarray} Leveraging the series element model in \eqref{eqngroup:genericSeriesForm}, we arrive at the multidimensional Ohm's law
\begin{IEEEeqnarray}{rCl}
\bmat{\mb{i} \\\mb{i}_{\mr{S}} }  = \mb{Y}_{\mr{net}} \bmat{\mb{v} \\ \mb{v}_{\mr{S}}},  \label{eqn:multidimensionalohmtildes}
\end{IEEEeqnarray}
where $\mb{v}_{\mr{S}}$  is the voltage at the slack bus and likewise $\mb{i}_{\mr{S}}$ is the current injection at the slack bus. Typically, $\mb{v}_{\mr{S}}$ is  equal to the symmetrical voltage  $\{1, 1\phase{ -120^{\circ}}, 1 \phase{120^{\circ}}\}$.
The matrix $\mb{Y}_{\mr{net}}$ can be constructed by block matrices leveraging the series model of transmission lines, SVRs, and transformers given in \eqref{eqngroup:genericSeriesForm} as follows:
\begin{subequations}
\label{eqn:ynet}
\begin{IEEEeqnarray}{rCl}
\mb{Y}_{\mr{net}} \left(\mc{J}_n , \mc{J}_n \right)  &=&  \sum\limits_{m \in \mc{N}_n} \mb{Y}_{nm}^{(n)}, \quad n \in \mc{N}   \label{eqn:ynetnn} \IEEEeqnarraynumspace\\
\mb{Y}_{\mr{net}} ( \mc{J}_n,\mc{J}_m) &=& - \mb{Y}_{nm}^{(m)},  \quad m \in \mc{N}_n. \label{eqn:ynetnm}
\end{IEEEeqnarray}
\end{subequations}
 Partitioning $\mb{Y}_{\mr{net}}$ yields
\begin{IEEEeqnarray}{rCl}
\bmat{\mb{i} \\ \mb{i}_{\mr{S}}} = \bmat{ \mb{Y} & \mb{Y}_{\mr{NS}} \\ \mb{Y}_{\mr{SN}}  & \mb{Y}_{\mr{SS}} } \bmat{ \mb{v} \\ \mb{v}_{\mr{S}}} \label{eqn:multiohm}
\end{IEEEeqnarray}
where $\mb{Y}$ is given by
\begin{subequations}
	\label{eqngroup:Y}
	\eq{rCl}{\mb{Y}(\mc{J}_n, \mc{J}_n) &=& \sum\limits_{m \in \mc{N}_n} \mb{Y}_{nm}^{(n)}, \quad n \in \mc{N} \setminus \{\mr{S}\} \label{eqn:Ynn} \\
		\mb{Y}(\mc{J}_n, \mc{J}_m) &=& - \mb{Y}_{nm}^{(m)}, \quad n \in \mc{N} \setminus \{\mr{S}\}, m \in \mc{N}_{n} \setminus \{\mr{S}\} \label{eqn:Ynm}  \IEEEeqnarraynumspace \\
		\mb{Y}(\mc{J}_n, \mc{J}_m) &=& \mb{O},  \quad n \in \mc{N} \setminus \{\mr{S}\}, m \notin \mc{N}_{n} \cup \{\mr{S}\}. \label{eqn:Ynmzero}}	
\end{subequations}
The matrix $\mb{Y}$ defined in \eqref{eqngroup:Y} is called the \emph{bus admittance matrix} of the network, and is also referred to as the Y-Bus.  Matrix $\mb{Y}$ is used for power flow analysis as we explain next.
	
From \eqref{eqn:multiohm}, it follows that $\mb{i}(\mb{v})  = \mb{Y} \mb{v} + \mb{Y}_{\mr{NS}} \mb{v}_{\mr{S}}$, 
where the dependence of $\mb{i}$ on $\mb{v}$ is made explicit. Due to the ZIP load model we can decompose $\mb{i}$ into three parts as follows:
\begin{IEEEeqnarray}{rCl}
\mb{i}_{\mr{PQ}}(\mb{v}) + \mb{i}_{\mr{I}} (\mb{v}) + \mb{i}_{\mr{Z}}(\mb{v}) = \mb{Y}\mb{v} + \mb{Y}_{\mr{NS}} \mb{v}_{\mr{S}}. \label{eqn:ipqizv}
\end{IEEEeqnarray}
From \eqref{eqn:YLn}, the constant-impedance currents is $\mb{i}_{\mr{Z}}(\mb{v}) = - \mb{Y}_{\mr{L}} \mb{v}$,  where $\mb{Y}_{\mr{L}} \in \mc{C}^{J \times J}$ is a block diagonal matrix with entries
\begin{IEEEeqnarray}{rCl}
\label{eqn:YLblk}
\mb{Y}_{\mr{L}} (\mc{J}_n, \mc{J}_n) = \mb{Y}_{\mr{L}_n}, \quad n \in \mc{N}\setminus \{\mr{S}\}. \label{eqn:YL} 
\end{IEEEeqnarray}

Using  $\mb{i}_{\mr{Z}}(\mb{v}) = - \mb{Y}_{\mr{L}} \mb{v}$ in \eqref{eqn:ipqizv} and rearranging yields
\begin{IEEEeqnarray}{rCl}
\mb{i}_{\mr{PQ}}(\mb{v}) + \mb{i}_{\mr{I}}(\mb{v})= (\mb{Y}+\mb{Y}_{\mr{L}}) \mb{v} + \mb{Y}_{\mr{NS}}\mb{v}_{\mr{S}}. \label{eqn:ipqiv}
\end{IEEEeqnarray}
If $\mb{Y}+\mb{Y}_{\mr{L}}$ is invertible, \eqref{eqn:ipqiv}  yields  a fixed-point equation for voltages, rendering the standard Z-Bus method as follows
\begin{IEEEeqnarray}{rCl}
\mb{v} [t+1]= \mb{Z} \left[  \mb{i}_{\mr{PQ}}(\mb{v}[t]) + \mb{i}_{\mr{I}} (\mb{v}[t]) \right] +\mb{w}, \label{eqn:zbus}
\end{IEEEeqnarray}
where $\mb{Z}= (\mb{Y} + \mb{Y}_{\mr{L}})^{-1}$, $\mb{w} = -\mb{Z}\mb{Y}_{\mr{NS}} \mb{v}_{\mr{S}}$, and $t$ is the iteration index. The Z-Bus method constitutes initializing $\mb{v}[0]$ typically to a flat voltage profile (i.e., $\mb{v}_n[0]:=\mb{v}_{\mr{S}}$) and running iteration \eqref{eqn:zbus}.  The solution in terms of voltages yields all power flows in the network.  The name of the method is derived from the inversion of $\mb{Y}+\mb{Y}_{\mr{L}}$ which yields $\mb{Z}$. In practice, $\mb{Y}+\mb{Y}_{\mr{L}}$ is not inverted, but its LU decomposition is computed which is unique if $\mb{Y}+\mb{Y}_{\mr{L}}$ is invertible.

It is well-known that $\mb{Y}$ is not always invertible due to transformer connections other than wye-g--wye-g. For all other transformer connections, matrices $\mb{Y}_2$, $\mb{Y}_3$, $\mb{Y}_4$, and $\mb{Y}_6$ in Table~\ref{table:transformerConnections} are rank-deficient, which may lead to a non-invertible $\mb{Y}$, depending on the position of the transformer in the feeder with respect to other elements.  The non-invertibility of $\mb{Y}$ can carry over to $\mb{Y}+\mb{Y}_{\mr{L}}$ as well. This is indeed the case for all the test feeders we investigate in Section \ref{sec:numtests}.

In order to numerically remedy the non-invertibility of $\mb{Y}$ for transformers other than wye-g--wye-g, \cite{Gorman1992} suggests adding a small shunt admittance from the isolated transformer sides to the ground.   Mathematically, this can be achieved by replacing $\mb{Y}_2$ and $\mb{Y}_4$ appearing in the self-admittances of Table~\ref{table:transformerConnections}  by  
\begin{subequations}
\label{eqngroup:Y2Y4prime} \eq{rCl}{\mb{Y}'_2&=& \mb{Y}_2+ \epsilon' \mb{I}, \label{eqn:Y2prime} \\
\mb{Y}'_4&=& \mb{Y}_4+ \epsilon'\mb{I}, \label{eqn:Y4prime} }
\end{subequations}
where $\epsilon'=\epsilon'_r-j\epsilon'_i$ with $\epsilon'_r>0$ is a small shunt admittance (compared to $|y_t|$).  If matrices $\mb{Y}_2$ and $\mb{Y}_4$  appear in the mutual admittances, then the modification is as follows: 
\begin{subequations}
	\label{eqngroup:Y2Y4doubleprime}
	\eq{rCl}{\mb{Y}''_2&=& \mb{Y}_2+ \epsilon'' \mb{I}, \label{eqn:Y2doubleprime} \\
		\mb{Y}''_4&=& \mb{Y}_4+ \epsilon''\mb{I}, \label{eqn:Y4doubleprime} }
\end{subequations}
where $\epsilon''=\epsilon''_r-j\epsilon''_i$ with $0 \le \epsilon''_r < \epsilon'_r$.
Doing so yields $\hat{\mb{Y}}$, a slightly modified version of $\mb{Y}$. \footnote{This modification is only intended to improve the numerical properties of the matrix $\mb{Y}$ in order to obtain steady-state power flow solutions. This is not to be confused with standard earthing practices through connection of a grounding  transformer such as wye-g--delta or zig-zag \cite{sallam2011electric}.}  In the next section, we prove the invertibility of the matrices $\hat{\mb{Y}}$ and $\hat{\mb{Y}}+\mb{Y}_{\mr{L}}$.  The  voltages can then be computed by using $\mb{Z}=(\hat{\mb{Y}}+\mb{Y}_{\mr{L}})^{-1}$  in \eqref{eqn:zbus}.  

A similar approach of modifying the transformer nodal admittances to handle the delta-side invertibility problem is employed in OpenDSS, where by default, a high reactance is connected from the delta terminal to the ground \cite{openDSSManual, openDSS,technoteTransformers}.

\section{Invertibility of  the bus admittance matrix $\mb{Y}$}
\label{sec:invertibility}
This section proves the invertibility of $\mb{Y}$, where to avoid exhaustive notation, the proof is based first on the following assumptions. 
The notation $\mb{S} \succ \mb{O}$ ($\mb{S} \succeq \mb{O}$) means matrix $\mb{S}$ is symmetric and positive (semi-)definite. 
 \renewcommand{\theenumi}{A\arabic{enumi}}
 \renewcommand{\labelenumi}{\theenumi)}
 \begin{enumerate}
 		\item \label{assumption:linecategories} The series elements belong to the following seven categories (nodal admittances are summarized in Table~\ref{table:sixcategories}): 
 		\renewcommand{\theenumii}{(\alph{enumii})}
 		\renewcommand{\labelenumii}{\theenumii}
 			\renewcommand{\theenumiii}{(\roman{enumiii})}
 			\renewcommand{\labelenumiii}{\theenumiii}
 			\makeatletter
 			\renewcommand\p@enumii{\theenumi}
 			\renewcommand\p@enumiii{\theenumi\theenumii}
 			\makeatother
 		\begin{enumerate}
 			\item \label{assumption:transmissionlines} The set $\mc{E}_{\mr{TL}}$ of transmission lines  with  nodal admittances given by \eqref{eqngroup:transmissionLineYtildes}. Two realistic properties verified in practical distribution feeders are assumed:
 			\begin{enumerate} 
 				\item \label{assumption:trsymmetric} Matrices $\mb{Z}_{nm}$ and $\mb{Y}_{nm}^{\mr{s}}$ are symmetric.
 				\item \label{assumption:trpd}It holds that $\mr{Re}[\mb{Z}_{nm}^{-1}] \succ \mb{O}$ and $\mr{Re}[\mb{Y}_{nm}^\mr{s}] \succeq \mb{O}$.
 			\end{enumerate} 
 			\item The set $\mc{E}_1$ of wye-g--wye-g transformers with nodal admittances given in Table~\ref{table:transformerConnections}. 
 			\item The set  $\mc{E}_2$ of wye-g--wye, wye--wye, and  delta--delta transformers, where in self-admittances the matrix $\mb{Y}'_2$ in~\eqref{eqn:Y2prime} replaces $\mb{Y}_2$ while in mutual admittances, the matrix $\mb{Y}''_2$ in~\eqref{eqn:Y2doubleprime}  replaces  $\mb{Y}_2$. \label{assumption:Y2}
 			\item The set $\mc{E}_{\mr{Yg\Delta}}$ of wye-g--delta transformers  where $\mb{Y}'_2$ in~\eqref{eqn:Y2prime} is used in place of $\mb{Y}_2$. \label{assumption:YgDelta} 
 			\item The set $\mc{E}_{\mr{Y\Delta}}$ of wye--delta transformers   where $\mb{Y}'_2$ in~\eqref{eqn:Y2prime} is used in place of $\mb{Y}_2$.  \label{assumtpion:Ydelta}
 			\item The set $\mc{E}_{\mr{O\Delta}}$ of open-delta--open-delta transformers  where in self-admittances, the matrix $\mb{Y}'_4$ in~\eqref{eqn:Y4prime} replaces $\mb{Y}_4$ while in mutual admittances, the matrix $\mb{Y}''_4$ in~\eqref{eqn:Y4doubleprime}  replaces  $\mb{Y}_4$. \label{assumption:Y4} 
 			 \item The set $\mc{E}_{\mr{R}}$  of  ideal  SVRs with nodal admittances given in \eqref{eqngroup:SVRIdealYtildes} and gain matrices $\mb{A}_{\mr{v}}$ and $\mb{A}_{\mr{i}}$ in Table~\ref{table:vrmatrices}. Notice that $\mb{A}_{\mr{v}}^{-1}=\mb{A}_{\mr{i}} ^T $ holds, and \ref{assumption:transmissionlines} is assumed for line $(n',m)$.
 			\begin{table}[t]
 			\centering
 				\begin{threeparttable}[b]
 				\renewcommand{\arraystretch}{1.3} 
 				\caption{Admittance matrices of series elements\tnote{$\dagger$}}
 				\small{
 					\begin{tabular}{|c|c|c|c|c|}
 						\hline
 						Edges & $\mb{Y}_{nm}^{(n)}$  & $\mb{Y}_{nm}^{(m)}$  & $\mb{Y}_{mn}^{(m)}$ & $\mb{Y}_{mn}^{(n)}$  \\ 
 						\hline
 						\multirow{2}{*}{$\mc{E}_{\mr{TL}}$} & $\frac{1}{2} \mb{Y}_{nm}^{\mr{s}}$  & \multirow{2}{*}{$\mb{Z}_{nm}^{-1}$} &  $\frac{1}{2} \mb{Y}_{nm}^{\mr{s}}$  & \multirow{2}{*}{$\mb{Z}_{nm}^{-1}$} \\
 						& $+\mb{Z}_{nm}^{-1}$ & & $+\mb{Z}_{nm}^{-1}$ &\\
 						\hline
 						$\mc{E}_1$ & $\mb{Y}_1$ & $\mb{Y}_1$ & $\mb{Y}_1$ & $\mb{Y}_1$ \\
 						\hline
 						$\mc{E}_2$ &  $\mb{Y}_2+\epsilon' \mb{I}$ &  $\mb{Y}_2+\epsilon'' \mb{I}$ &  $\mb{Y}_2+\epsilon' \mb{I}$ &  $\mb{Y}_2+\epsilon'' \mb{I}$\\
 						\hline 
 						$\mc{E}_{\mr{Yg\Delta}}$ & $\mb{Y}_1$ & $-\mb{Y}_3$ &  $\mb{Y}_2+\epsilon' \mb{I}$ & $-\mb{Y}_3^T$ \\
 						\hline
 						$\mc{E}_{\mr{Y\Delta}}$ & $\mb{Y}_2+\epsilon' \mb{I}$ & $-\mb{Y}_3$ & $\mb{Y}_2+\epsilon' \mb{I}$ & $-\mb{Y}_3^T$ \\
 						\hline 
 						$\mc{E}_{\mr{O\Delta}}$ & $\mb{Y}_4 + \epsilon' \mb{I}$ & $\mb{Y}_4 + \epsilon'' \mb{I}$ & $\mb{Y}_4 + \epsilon' \mb{I}$ & $\mb{Y}_4 + \epsilon'' \mb{I}$ \\
 						\hline
 						\multirow{2}{*}{$\mc{E}_{\mr{R}}$} & $\mb{A}_{\mr{i}}  \bigl(\frac{1}{2} \mb{Y}_{n'm}^{\mr{s}} \bigr.$ & \multirow{2}{*}{$\mb{A}_{\mr{i}} \mb{Z}_{n'm}^{-1}$} & $\frac{1}{2} \mb{Y}_{n'm}^{\mr{s}}$ & \multirow{2}{*}{$\mb{Z}_{n'm}^{-1}\mb{A}_{\mr{i}}^T$} \\
 						&  $\bigl.+\mb{Z}_{n'm}^{-1}\bigr)\mb{A}_{\mr{i}}^T$ & & $+ \mb{Z}_{n'm}^{-1}$ & \\
 						\hline  						$\mc{E}_{\mr{OY\Delta}}$ & $\mb{Y}_5 $ & $-\mb{Y}_6 $ & $\mb{Y}_4 + \epsilon' \mb{I}$ & $-\mb{Y}_6^T$ \\
 						\hline
 					\end{tabular}}
 					\label{table:sixcategories}
 				\end{threeparttable}
 				\begin{tablenotes}
 		\item $^\dagger$ In case of missing phases, admittance matrices select $\Omega_{nm}$, cf.~\eqref{eqngroup:transmissionLineYtildes}.
 				\end{tablenotes}
 				\end{table}
  	\end{enumerate}
Notice that $\mc{E}= \mc{E}_{\mr{TL}} \cup \mc{E}_1 \cup \mc{E}_2 \cup \mc{E}_{\mr{Yg\Delta}} \cup \mc{E}_{\mr{Y\Delta}} \cup \mc{E}_{\mr{O\Delta}} \cup \mc{E}_{\mr{R}}$.
  		\item \label{assumption:connectedness}The graph $(\mc{N},\mc{E})$ is connected.
 	\item \label{assumption:nonidealtransformers} For all transformers it holds that $\mr{Re}[y_t] >0$, i.e., the transformers are non-ideal.  
 	\item \label{assumption:threephases} Nodes and edges have three available phases, that is $\Omega_n=\Omega_{nm}=\Omega_m=\{a,b,c\}$ for all $n \in \mc{N}$ and $(n,m) \in \mc{E}$.  The extension to the case that includes transmission lines with missing phases and open-wye--open-delta transformers is provided afterwards.
 \end{enumerate}
 \renewcommand{\labelenumi}{\theenumi)}
A few comments are in order here regarding Assumption \ref{assumption:transmissionlines}. The symmetricity of  $\mb{Z}_{nm}$ and $\mb{Y}_{nm}^{\mr{s}}$ holds from electromagnetic considerations, that is, the coupling between phases $a$ and $b$ is the same as the coupling between $b$ and $a$.  Interestingly, we verified that matrix $\mr{Re}[\mb{Z}_{nm}]$ is typically diagonally dominant with positive diagonal entries for all distribution feeders. Hence, it holds that $\mr{Re}[\mb{Z}_{nm}] \succ \mb{O}$, and due to Lemma \ref{lemma:usefulalgebra} in Appendix~\ref{sec:appendixLemma1}, $\mb{Z}_{nm}^{-1}$ exists, and  $\mr{Re}[\mb{Z}_{nm}^{-1}] \succ \mb{O}$.  Assumption \ref{assumption:trpd} is though weaker than diagonal dominance. It is worth emphasizing that \ref{assumption:transmissionlines} allows for \emph{untransposed} lines.

The invertibility proof relies on Lemma~\ref{lemma:usefulalgebra} in Appendix~\ref{sec:appendixLemma1}. To use this lemma, it is first proved that $\mb{Y}$ is symmetric (Theorem~\ref{lemma:Ysymmetric}) and that    $\mr{Re}[\mb{Y}] \succ \mb{O}$ (Theorem~\ref{lemma:realYpositivedefinite}). 

 \begin{theorem} \label{lemma:Ysymmetric} Under \ref{assumption:linecategories},  $\mb{Y}$ defined in \eqref{eqngroup:Y} is symmetric. 
 \end{theorem}
\begin{IEEEproof}
It suffices to prove that $\mb{Y}_{nm}^{(n)} = (\mb{Y}_{nm}^{(n)})^T$ and $\mb{Y}_{nm}^{(m)} = (\mb{Y}_{mn}^{(n)})^T$ for all $n \in \mc{N}\setminus\{\mr{S}\}$ and $m \in \mc{N}_n \setminus \{\mr{S}\}$ [cf.~\eqref{eqngroup:Y}] where the admittance matrices are given in Table~\ref{table:sixcategories}. The relationships hold for $(n,m) \in \mc{E}_{\mr{TL}}$ due to  \ref{assumption:trsymmetric}; for $(n,m) \in  \mc{E}_1 \cup \mc{E}_2  \cup \mc{E}_{\mr{Yg\Delta}} \cup \mc{E}_{\mr{Y\Delta}} \cup \mc{E}_{\mr{O\Delta}}$ due to~\eqref{eqngroup:y1y2y3}; and for  $(n,m) \in \mc{E}_{\mr{R}}$ due to~\ref{assumption:trsymmetric} and using $\mb{A}_{\mr{v}}^{-1}=\mb{A}_{\mr{i}} ^T $. 
\end{IEEEproof} 
Let $\mc{E}_{\mr{S}}$ collect the edges connected to the slack bus.   Such an edge could be a transmission line, a voltage regulator, or a three-phase transformer according to Table~\ref{table:sixcategories}. 
 \begin{theorem} \label{lemma:realYpositivedefinite} Under \ref{assumption:linecategories}--\ref{assumption:threephases}, it holds that  $\mr{Re}[\mb{Y}] \succ \mb{O}$. 
 \end{theorem}
 \begin{IEEEproof} {\normalsize The goal is to prove that $\mb{x}^T \mr{Re} [\mb{Y}]\mb{x} >0$ for any nonzero real vector $\mb{x} \in \mbb{R}^J$. Since $\mb{x}$ is real, notice that $\mb{x}^T \mr{Re} [\mb{Y}]\mb{x}  = \mr{Re}[ \mb{x}^T \mb{Y} \mb{x} ]$.  By the definition of $\mb{Y}$ in \eqref{eqngroup:Y}, and defining $\mb{x}_n$ as a vector that selects indices of $\mb{x}$ corresponding to $\mc{J}_n$, $\mb{x}^T  \mr{Re}[\mb{Y}] \mb{x}$ is written as follows:
 
 ~}
 \vspace{-1em}
 		{\small \eq{rCl}{\IEEEeqnarraymulticol{3}{l}{ \mb{x}^T \mr{Re}[\mb{Y}] \mb{x}=  \mr{Re}\left[\sum\limits_{n \in  \mc{N}\setminus \{\mr{S}\} }  \hspace{-0.2cm}\mb{x}_n^T \left(\sum \limits_{m \in \mc{N}_n} \hspace{-0.2cm} \mb{Y}_{nm}^{(n)}  \mb{x}_n  - \hspace{-0.4cm}\sum \limits_{m \in \mc{N}_n\setminus \{\mr{S}\}} \hspace{-0.2cm} \mb{Y}_{nm}^{(m)}  \mb{x}_m\right) \right]}\notag \\
 		&=&  \sum\limits_{n \in  \mc{N}_{\mr{S}}} \hspace{-0.1cm}  \mb{x}_n^T \mr{Re}[\mb{Y}_{n\mr{S}}^{(n)}]  \mb{x}_n +\hspace{-0.5cm} \sum\limits_{\substack{n \in \mc{N} \setminus \{\mr{S}\} \\  m \in \mc{N}_n \setminus \{\mr{S}\} }}\hspace{-0.4cm} \mr{Re} \left[ \mb{x}_n^T \mb{Y}_{nm}^{(n)}  \mb{x}_n  - \mb{x}_n^T\mb{Y}_{nm}^{(m)}  \mb{x}_m\right]\IEEEeqnarraynumspace \notag \\ 
 		&=&  \sum\limits_{n \in  \mc{N}_{\mr{S}}} \hspace{-0.1cm} \mb{x}_n^T \mr{Re}[\mb{Y}_{n\mr{S}}^{(n)}]  \mb{x}_n +\hspace{-0.5cm} \sum\limits_{(n,m) \in \mc{E} \setminus \mc{E}_{\mr{S}}} \hspace{-0.5cm}\mr{Re}\left[\mb{x}_n^T \mb{Y}_{nm}^{(n)}  \mb{x}_n  - \mb{x}_n^T\mb{Y}_{nm}^{(m)}  \mb{x}_m \right. \notag \\
 		&&  \hspace{27mm} \left. + \: \mb{x}_m ^T  \mb{Y}_{mn}^{(m)}\mb{x}_m - \mb{x}_m^T \mb{Y}_{mn}^{(n)}\mb{x}_n\right]. \IEEEeqnarraynumspace \label{eqn:xtyx}} }
In \eqref{eqn:xtyx}, replace the second summation using Lemmas~\ref{lemma:utWv} and~\ref{lemma:utWXZv} in Appendix~\ref{sec:appendixLemma1} as follows. For  $\mc{E}_{\mr{TL}} \setminus \mc{E}_{\mr{S}}$, apply Lemma \ref{lemma:utWv} with $\mb{u}= \mb{x}_n$, $\mb{v}= \mb{x}_m$, and $\mb{W}= \mr{Re}[\mb{Z}_{nm}^{-1}]$.
 For $ \mc{E}_1 \setminus \mc{E}_{\mr{S}}$, apply Lemma \ref{lemma:utWv} with $\mb{u}= \mb{x}_n$, $\mb{v}= \mb{x}_m$, $\mb{W}= \mr{Re}[\mb{Y}_1]$.
For $ \mc{E}_2 \setminus \mc{E}_{\mr{S}}$, apply Lemma \ref{lemma:utWv} with $\mb{u}= \mb{x}_n$, $\mb{v}= \mb{x}_m$, and $\mb{W}= \mr{Re}[\mb{Y}_2]+ \epsilon''_r \mb{I}$.
 For $ \mc{E}_{\mr{Yg\Delta}} \setminus \mc{E}_{\mr{S}}$  apply Lemma \ref{lemma:utWXZv} with $\mb{u}= \mb{x}_n$, $\mb{v}= \mb{x}_m$, and $\mb{W}= \mr{Re}[\mb{Y}_1]$, $\mb{X}= -\mr{Re}[\mb{Y}_3]$, $\mb{Z}= \mr{Re}[\mb{Y}_2]$.
 For $\mc{E}_{\mr{Y\Delta}} \setminus \mc{E}_{\mr{S}}$, apply Lemma \ref{lemma:utWXZv} with $\mb{u}= \mb{x}_n$, $\mb{v}= \mb{x}_m$, and $\mb{W}=\mb{Z}= \mr{Re}[\mb{Y}_2]$, $\mb{X}= -\mr{Re}[\mb{Y}_3]$.
For $ \mc{E}_4 \setminus \mc{E}_{\mr{S}}$, apply Lemma \ref{lemma:utWv} with $\mb{u}= \mb{x}_n$, $\mb{v}= \mb{x}_m$, and $\mb{W}= \mr{Re}[\mb{Y}_4]+ \epsilon''_r \mb{I}$.
 For $ \mc{E}_{\mr{R}} \setminus \mc{E}_{\mr{S}}$,  apply Lemma \ref{lemma:utWv} with $\mb{u}= \mb{A}_{\mr{i}}^T\mb{x}_n$, $\mb{v}= \mb{x}_m$, and $\mb{W}= \mr{Re}[\mb{Z}_{n'm}^{-1}]$. After the previous replacements, we obtain
 
{\small \begin{IEEEeqnarray}{rCl}
	 	\IEEEeqnarraymulticol{3}{l}{\mb{x}^T \mr{Re}[\mb{Y}] \mb{x}=  \sum\limits_{n \in  \mc{N}_{\mr{S}}}  \mb{x}_n^T \mr{Re}[\mb{Y}_{n\mr{S}}^{(n)}]  \mb{x}_n}  \notag \\
	\quad  &+& \:\hspace{-0.62cm} \sum \limits_{(n,m) \in \mc{E}_{\mr{TL}} \setminus \mc{E}_{\mr{S}}} \hspace{-0.58cm}(\mb{x}_n- \mb{x}_m)^T\mr{Re}[\mb{Z}_{nm}^{-1}] (\mb {x}_n - \mb{x}_m)   \IEEEeqnarraynumspace \notag \\
	 && + \: \frac{1}{2}\mb{x}_n^T \mr{Re}[\mb{Y}_{nm}^{\mr{s}}] \mb{x}_n^T + \frac{1}{2}\mb{x}_m^T \mr{Re}[\mb{Y}_{nm}^{\mr{s}}] \mb{x}_m^T  \IEEEeqnarraynumspace \notag \\
	\quad  &+& \: \hspace{-0.62cm}\sum \limits_{(n,m) \in \mc{E}_1 \setminus \mc{E}_{\mr{S}}}  \hspace{-0.50cm}(\mb{x}_n- \mb{x}_m)^T \mr{Re}[\mb{Y}_1] (\mb {x}_n - \mb{x}_m)  \IEEEeqnarraynumspace \notag \\
	 \quad &+& \: \hspace{-0.62cm} \sum \limits_{(n,m) \in \mc{E}_2 \setminus \mc{E}_{\mr{S}}} \hspace{-0.50cm}(\mb{x}_n- \mb{x}_m)^T (\mr{Re}[\mb{Y}_2]+\epsilon_r'' \mb{I}) (\mb {x}_n - \mb{x}_m) \notag \\
	 	&& + \: (\epsilon'_r-\epsilon''_r)(\mb{x}_n^T\mb{x}_n+\mb{x}_m^T\mb{x}_m)  \IEEEeqnarraynumspace   \notag \\
		\quad &+& \:  \hspace{-0.62cm}\sum \limits_{(n,m) \in \mc{E}_{\mr{Yg\Delta}}\setminus \mc{E}_{\mr{S}}}\hspace{-0.6cm} \bmat{\mb{x}_n^T & \mb{x}_m^T}  	\mb{G}_{\mr{Yg\Delta}_{nm}} \bmat{\mb{x}_n \\ \mb{x}_m} +   \epsilon'_r \mb{x}_m^T \mb{x}_m  \IEEEeqnarraynumspace \notag \\
		\quad &+& \: \hspace{-0.62cm} \sum \limits_{(n,m) \in \mc{E}_{\mr{Y\Delta}} \setminus \mc{E}_{\mr{S}}} \hspace{-0.6cm} \bmat{\mb{x}_n^T & \mb{x}_m^T}  	\mb{G}_{\mr{Y\Delta}_{nm}} \bmat{\mb{x}_n \\ \mb{x}_m} \notag \\
			 && + \:  \epsilon'_r\mb{x}_n^T \mb{x}_n + \epsilon'_r \mb{x}_m^T  \mb{x}_m \IEEEeqnarraynumspace \notag \\
			 		 	 \quad &+& \: \hspace{-0.62cm} \sum \limits_{(n,m) \in \mc{E}_4 \setminus \mc{E}_{\mr{S}}} \hspace{-0.50cm}(\mb{x}_n- \mb{x}_m)^T (\mr{Re}[\mb{Y}_4]+\epsilon_r'' \mb{I}) (\mb {x}_n - \mb{x}_m) \notag \\
		&& + \: (\epsilon'_r-\epsilon''_r)(\mb{x}_n^T\mb{x}_n+\mb{x}_m^T\mb{x}_m) \IEEEeqnarraynumspace  \notag \\
	&+& \: \hspace{-0.62cm} \sum \limits_{(n,m) \in \mc{E}_{\mr{R}} \setminus \mc{E}_{\mr{S}}} \hspace{-0.6cm}\Bigl\{ (\mb{A}_{\mr{i}}^T\mb{x}_n- \mb{x}_m)^T \mr{Re}[\mb{Z}_{n'm}^{-1}] (\mb{A}_{\mr{i}}^T\mb {x}_n - \mb{x}_m) \Bigr.  \IEEEeqnarraynumspace \notag  \\
	 && \Bigl. + \: \frac{1}{2}\mb{x}_n^T \mb{A}_{\mr{i}}\mr{Re}[\mb{Y}_{n'm}^{\mr{s}}] \mb{A}_{\mr{i}}^T \mb{x}_n +\frac{1}{2} \mb{x}_m^T \mr{Re}[\mb{Y}_{n'm}^{\mr{s}}] \mb{x}_m\Bigr\}
	\label{eqn:xtGxbreakdown}
\end{IEEEeqnarray}}
{\normalsize where we set  $\mb{G}_{\mr{Yg\Delta}_{nm}} = \smat{\mr{Re}[\mb{Y}_1] & \mr{Re}[\mb{Y}_3] \\ \mr{Re}[\mb{Y}_3^T] & \mr{Re}[\mb{Y}_2]}$, and  $\mb{G}_{\mr{Y\Delta}_{nm}} = \smat{\mr{Re}[\mb{Y}_2] & \mr{Re}[\mb{Y}_3] \\ \mr{Re}[\mb{Y}_3^T] & \mr{Re}[\mb{Y}_2]}$}. 
It holds that $\mb{G}_{\mr{Yg\Delta}_{nm}}, 	\mb{G}_{\mr{Y\Delta}_{nm}} \succeq \mb{O}$ since their nonzero eigenvalues are respectively  $\{\mr{Re}[y_t], 2\mr{Re}[y_t]\}$, and $\{2\mr{Re}[y_t]\}$; and $\mr{Re}[y_t] > 0$ for all transformers.  Furthermore,  $\mr{Re}[\mb{Z}_{nm}^{-1}] \succ \mb{O}$ for transmission lines, and it holds that $\mr{Re}[\mb{Y}_1] \succ \mb{O}$, $\mr{Re}[\mb{Y}_2] \succeq \mb{O}$, $\mr{Re}[\mb{Y}_4] \succeq \mb{O}$, and $\epsilon_r'>\epsilon_r''$.  Moreover, $\mb{Y}_{n\mr{S}}^{(n)}$ is the self-admittance of any of the edges in Table~\ref{table:sixcategories}, and therefore,  $\mr{Re}[\mb{Y}_{n\mr{S}}^{(n)}] \succ \mb{O}$ holds.  The expression in \eqref{eqn:xtGxbreakdown} is thus nonnegative for any real vector $\mb{x}$. This proves that $\mr{Re}[\mb{Y}] \succeq \mb{O}$. At this point, to prove positive definiteness, we show that whenever $\mb{x}^T \mr{Re}[\mb{Y}] \mb{x}=0$ holds, then $\mb{x}=\mb{0}$.  In particular, if $\mb{x}^T \mr{Re}[\mb{Y}] \mb{x}=0$, then every quadratic form in \eqref{eqn:xtGxbreakdown} must be zero and we have the following implications:
 \renewcommand{\theenumi}{R\arabic{enumi}}
  \renewcommand{\labelenumi}{(\theenumi)}
\begin{enumerate}[leftmargin=3\parindent]
	\item \label{enum:resultnomissing1}  $\mr{Re}[\mb{Y}_{n\mr{S}}^{(n)}] \succ \mb{O} \Rightarrow \mb{x}_n=\mb{0}$, $n \in \mc{N}_{\mr{S}}$.
	\item \label{enum:resultnomissing2} $\mr{Re}[\mb{Z}_{nm}^{-1}] \succ \mb{O} \Rightarrow \mb{x}_n=\mb{x}_m$, $(n,m) \in \mc{E}_{\mr{TL}} \setminus \mc{E}_{\mr{S}}$.
	\item \label{enum:resultnomissingY1}  $\mr{Re}[\mb{Y}_1] \succ \mb{O}\Rightarrow \mb{x}_n=\mb{x}_m$,  $(n,m) \in \mc{E}_1 \setminus \mc{E}_{\mr{S}}$.
		\item \label{enum:resultnomissingY2}  $\epsilon_r'-\epsilon_r'' > 0 \Rightarrow \mb{x}_n=\mb{x}_m=\mb{0}$,  $(n,m) \in \mc{E}_2 \setminus \mc{E}_{\mr{S}}$. 
	\item  \label{enum:resultnomissingYgDelta}$\epsilon'_r >0 \Rightarrow \mb{x}_m=\mb{0}$; and $\mr{Re}[\mb{Y}_1] \succ \mb{O} \Rightarrow \mb{x}_n=\mb{0}$, $(n,m) \in \mc{E}_{\mr{Yg\Delta}} \setminus \mc{E}_{\mr{S}}$.
		\item \label{enum:resultnomissingYdelta} $\epsilon'_r >0 \Rightarrow \mb{x}_n=\mb{x}_m=\mb{0}$, $(n,m) \in \mc{E}_{\mr{Y\Delta}} \setminus \mc{E}_{\mr{S}}$.
		\item \label{enum:resultnomissingY4}  $\epsilon_r' - \epsilon_r'' >0 \Rightarrow \mb{x}_n=\mb{x}_m=\mb{0}$,  $(n,m) \in \mc{E}_4 \setminus \mc{E}_{\mr{S}}$.
	\item \label{enum:resultnomissingRegs} $\mr{Re}[\mb{Z}_{n'm}^{-1}] \succ \mb{O} \Rightarrow \mb{x}_m= \mb{A}_{\mr{i}} ^T\mb{x}_n$,  $(n,m) \in \mc{E}_{\mr{YVR}} \setminus \mc{E}_{\mr{S}}$.
\end{enumerate}
Since  the graph is connected,  $\mb{x}^T \mr{Re}[\mb{Y}] \mb{x}= 0$ implies that $\mb{x}_n=\mb{0}$ for all $n$. We conclude that $\mr{Re}[\mb{Y}] \succ \mb{O}$.  
 \end{IEEEproof}
\begin{remark}
The concluding part of the proof in Theorem \ref{lemma:realYpositivedefinite} was inspired by a similar argument in \cite{WangBernsteinBoudecPaolone2016} that proved the invertibility of the single-phase Y-Bus.
\end{remark}

\begin{theorem} \label{theorem:Yinv}
Under \ref{assumption:linecategories}--\ref{assumption:threephases}, matrix $\mb{Y}$ is invertible.
\end{theorem}
\begin{IEEEproof}
Combine Theorems \ref{lemma:Ysymmetric} and~\ref{lemma:realYpositivedefinite} with Lemma \ref{lemma:usefulalgebra}. 
\end{IEEEproof}
This section proved  the invertibility of  the Y-Bus under Assumptions \ref{assumption:linecategories}--\ref{assumption:threephases}. It is worth emphasizing that the  proof is modular, that is, the invertibility of the Y-Bus is shown to be dependent only on the  properties of block admittances of the edges.  In other words, regardless of the number and order of elements, and of the topology (radial or meshed), the Y-Bus is invertible, as long as the series elements adhere to Assumptions \ref{assumption:linecategories}--\ref{assumption:threephases}.  The next section relaxes some of these assumptions.

\section{Extensions}
\label{sec:extension}
In this section, extensions to Theorem \ref{theorem:Yinv} that guarantee invertibility of the Y-Bus  under more practical considerations are provided. Specifically, Section \ref{sec:invertyYmissingphase} covers typical distribution feeders that consist of  one-, two-, and three-phase laterals and open-wye--open-delta transformers. The Z-Bus method of~\eqref{eqn:zbus} relies on the invertibility of $\mb{Y}+\mb{Y}_{\mr{L}}$, which is asserted in Section \ref{sec:invertyyl}. Section \ref{sec:invertimag} provides a proof of the Y-Bus invertibility using only the imaginary parts of $\epsilon'$ and $\epsilon''$; thereby allowing the inclusion of transformers with zero ohmic losses. This is indeed useful for simplistic models  where transformers are only represented by their series reactances. Finally, Section~\ref{sec:altxfm} considers the Y-Bus invertibility under alternative practical transformer connections.

\subsection{Invertibility under transmission lines with missing phases and open-wye--open-delta transformers}
\label{sec:invertyYmissingphase}
The ensuing Theorem~\ref{theorem:missingphases} handles circuits that include transmission lines with missing phases and open-wye--open-delta transformers. To this end, Assumptions ~\ref{assumption:linecategories} and ~\ref{assumption:threephases} are modified as follows:
\renewcommand{\theenumi}{A\arabic{enumi}$'$}
\renewcommand{\labelenumi}{\theenumi)}
\begin{enumerate}
\item The series elements in Assumption \eqref{assumption:linecategories} could also include the set $\mc{E}_{\mr{OY\Delta}}$ as the set  of open-wye--open-delta transformers where $\mb{Y}_4'$ in~\eqref{eqn:Y4prime} is used in place of $\mb{Y}_4$.  \label{assumption:open-wye--open-delta}  
\setcounter{enumi}{3}
	\item \label{assumption:missingphases} For any node $n \in \mc{N}$ and phase $\phi \in \Omega_n$, a path exists from  $\mr{S}$ to $n$   where edges   have $\phi$ in their phase set. 
\end{enumerate}
Typical distribution systems comprise a main three-phase feeder with one-, two-, or three-phase laterals, and therefore, \ref{assumption:missingphases}  is satisfied in practice.
\begin{theorem} \label{theorem:missingphases}
	The conclusion of Theorem \ref{lemma:realYpositivedefinite} holds under \ref{assumption:open-wye--open-delta}, \ref{assumption:connectedness}, \ref{assumption:nonidealtransformers}, and \ref{assumption:missingphases}.
\end{theorem}
\begin{IEEEproof}
The only modification will occur in the term \eqref{eqn:xtGxbreakdown} where for transmission lines $(n,m) \in \mc{E}_{\mr{TL}}$ and   SVRs $(n,m) \in \mc{E}_{\mr{R}}$,  we	substitute $\mb{x}_n$ and $\mb{x}_m$ respectively with $\mb{x}_n(\Omega_{nm})$ and $\mb{x}_m(\Omega_{nm})$.  For open-wye--open-delta transformers$(n,m) \in \mc{E}_{\mr{OY\Delta}}$ we substitute $\mb{x}_n$  with $\mb{x}_n(\Omega_{nm})$.  Moreover, we include the following term for  open-wye--open-delta transformers$(n,m) \in \mc{E}_{\mr{OY\Delta}}$ in \eqref{eqn:xtGxbreakdown}
	
	\eq{rCl}{	\sum \limits_{(n,m) \in \mc{E}_{\mr{OY\Delta}} \setminus \mc{E}_{\mr{S}}} \hspace{-0.6cm} \bmat{\mb{x}_n(\Omega_{nm})^T & \mb{x}_m^T}  	\mb{G}_{\mr{OY\Delta}_{nm}} \bmat{\mb{x}_n(\Omega_{nm}) \\ \mb{x}_m} ,   \IEEEeqnarraynumspace \label{eqn:oydeltatr1}}

\noindent \normalsize where $\mb{G}_{\mr{OY\Delta}_{nm}}= \smat{ \mr{Re}[\mb{Y}_5] & \mr{Re}[\mb{Y}_6] \\ \mr{Re}[\mb{Y}_6^T] &  \mr{Re}[\mb{Y}_4]+\epsilon_r' \mb{I}}$.  Notice that $\mr{Re}[\mb{Y}_5] =\mr{Re}[y_t] \mb{I}_2$. Moreover, the non-zero eigen-values of $\mr{Re}[\mb{Y}_4]$ are $\{\frac{1}{3}\mr{Re}[y_t],\mr{Re}[y_t]\}$ and thus $\mr{Re}[\mb{Y}_4] \succeq \mb{O}$. 
Thus, the Schur complement of $\mb{G}_{\mr{OY\Delta}}$ is given by 
\eq{rCl}{\mr{Re}[\mb{Y}_4]+\epsilon_r' \mb{I}-\mr{Re}[\mb{Y}_6^T] \mr{Re}[\mb{Y}_5]^{-1} \mr{Re}[\mb{Y}_6].  \label{eqn:shurcomplementGYOD}}
Furthermore, the following identity holds: $\mr{Re}[\mb{Y}_6^T] \mr{Re}[\mb{Y}_6]=\mr{Re}[y_t]\mr{Re}[\mb{Y}_4]$. Therefore, if $\epsilon'_r >0$, the matrix $\mr{Re}[\mb{Y}_4]+\epsilon_r' \mb{I}$ as well as the Schur complement of $\mb{G}_{\mr{OY\Delta}}$  calculated in \eqref{eqn:shurcomplementGYOD} are positive definite.   This implies that $\mb{G}_{\mr{OY\Delta}}$ is positive definite if $\epsilon'_r >0$. 
Thus, we conclude that
 \renewcommand{\theenumi}{R\arabic{enumi}$'$}
 \renewcommand{\labelenumi}{(\theenumi)}
\begin{enumerate}[leftmargin=3\parindent]
\setcounter{enumi}{1}
\item \label{enum:missingLines}  For  $(n,m)\in \mc{E}_{\mr{TL}} \setminus \mc{E}_{\mr{S}}$, $\mb{x}_{n}(\Omega_{nm})= \mb{x}_m(\Omega_{nm})$ 
		\setcounter{enumi}{7}
\item \label{enum:missingRegs}  For  $(n,m) \in \mc{E}_{\mr{R}} \setminus \mc{E}_{\mr{S}}$, $\mb{x}_m(\Omega_{nm})=\mb{A}_{\mr{i}}^T\mb{x}_n(\Omega_{nm})$. 
\item \label{enum:missingOYdelta} For  $(n,m) \in \mc{E}_{\mr{OY\Delta}} \setminus \mc{E}_{\mr{S}}$, since $\mb{G}_{\mr{OY\Delta}} \succ \mb{O} \Rightarrow \mb{x}_m=\mb{0}$ and $\mb{x}_n(\Omega_{nm})=\mb{0}$. 
	\end{enumerate}
	\renewcommand{\theenumi}{\arabic{enumi}}
Consider an edge $(n,m) \in ( \mc{E}_{\mr{TL}} \cup \mc{E}_{\mr{R}} \cup \mc{E}_{{\mr{OY\Delta}}}) \setminus \mc{E}_{\mr{S}}$ with  $m \notin \mc{N}_{\mr{S}}$ and $\phi \in \Omega_m$.  Due to \ref{assumption:missingphases}, there exists a path from the slack bus $\mr{S}$  to node $m$ where all edges  must have $\phi$ in their phase set.  This path is denoted by $\mc{P}=\{\mr{S}, j, \ldots, m\}$ where $j \in \mc{N}_{\mr{S}}$. Based on \eqref{enum:resultnomissing1}, $\mb{x}_j(\{\phi\})=0$. Based on \eqref{enum:missingLines}, \eqref{enum:missingRegs},  \eqref{enum:missingOYdelta}, \eqref{enum:resultnomissingY1}--\eqref{enum:resultnomissingY4}, we conclude that $\mb{x}_m (\{\phi\})=0$.
	\end{IEEEproof}

\subsection{Invertibility of $\mb{Y}+\mb{Y}_{\mr{L}}$}
\label{sec:invertyyl}
Prior to presenting that $\mb{Y}+\mb{Y}_{\mr{L}}$ is invertible, the ensuing Theorem \ref{lemma:YLnsymmetric} is first presented that  brings out some useful properties of $\mb{Y}_{\mr{L}}$. 
\begin{theorem} \label{lemma:YLnsymmetric} Matrix $\mb{Y}_{\mr{L}}$ defined in \eqref{eqn:YL} is symmetric, and it holds that $\mr{Re}[\mb{Y}_{\mr{L}}] \succeq \mb{O}$. 
\end{theorem}
\begin{IEEEproof}
	Since $\mb{Y}_{\mr{L}}$ is  block diagonal [cf.~\eqref{eqn:YLblk}], it suffices to prove that $\mb{Y}_{\mr{L}_n}$ in \eqref{eqngroup:YLnspecific} is symmetric and that $\mr{Re}[\mb{Y}_{\mr{L}_n}] \succeq \mb{O}$ for $n \in \mc{N}\setminus \{\mr{S}\}$.  For wye loads, $\mb{Y}_{\mr{L}_n}$ is diagonal due to \eqref{eqn:YLnwye}; and  since  $\mr{Re}[y_{L_n}^{\phi}] \ge 0$, it holds that $\mr{Re}[\mb{Y}_{\mr{L}_n}] \succeq \mb{O}$.  For delta loads, symmetricity of $\mb{Y}_{\mr{L}_n}$ follows from \eqref{eqn:YLndeltaoffdiag}; also 
	$\mr{Re}[\mb{Y}_{\mr{L}_n}]$ is symmetric diagonally dominant  with nonnegative diagonal entries, and thus $\mr{Re}[\mb{Y}_{\mr{L}_n}] \succeq \mb{O}$. 
\end{IEEEproof}  

\begin{theorem} \label{theorem:YYLinvert}
The matrix $\mb{Y}+\mb{Y}_{\mr{L}}$ with $\mb{Y}$ defined in \eqref{eqngroup:Y} and $\mb{Y}_{\mr{L}}$ defined in \eqref{eqn:YL} is invertible under \ref{assumption:open-wye--open-delta}, \ref{assumption:nonidealtransformers}, \ref{assumption:connectedness}, \ref{assumption:missingphases}.
\end{theorem} 
\begin{IEEEproof}
Let $\mb{Y}+\mb{Y}_{\mr{L}}= \mb{G}+j \mb{B}$. Theorems~\ref{lemma:YLnsymmetric}, \ref{lemma:Ysymmetric}, and~\ref{lemma:realYpositivedefinite} or~\ref{theorem:missingphases} prove that $\mb{G} =\mr{Re}[\mb{Y}+\mb{Y}_{\mr{L}}]=\mr{Re}[\mb{Y}]+\mr{Re}[\mb{Y}_{\mr{L}}] \succ \mb{O}$ and $\mb{B} = \mr{Im}[\mb{Y}+\mb{Y}_{\mr{L}}]= \mr{Im}[\mb{Y}^T+\mb{Y}_{\mr{L}}^T]=\mb{B}^T$.  By Lemma~\ref{lemma:usefulalgebra}, matrix $\mb{Y}+\mb{Y}_{\mr{L}}$ is invertible.
\end{IEEEproof}

The previous analysis reveals that transformers with $\mb{Y}_2$ or $\mb{Y}_4$ as their self-admittance compromise the invertibility of $\mb{Y}+\mb{Y}_{\mr{L}}$. The modification of $\mb{Y}_2$ or similarly $\mb{Y}_4$  in \eqref{eqngroup:Y2Y4prime} and \eqref{eqngroup:Y2Y4doubleprime} guarantees invertibility in Theorem \ref{theorem:Yinv}.   It is alternatively possible to achieve invertibility of $\mb{Y}+\mb{Y}_{\mr{L}}$ when $\mb{Y}_{\mr{L}}$ automatically adds a constant-impedance load with positive definite real part to a side of the transformer whose self-admittance is $\mb{Y}_2$ or $\mb{Y}_4$. This observation is summarized in the following corollary. 
\begin{corollary}
Matrix $\mb{Y}+\mb{Y}_{\mr{L}}$  is invertible if  matrices $\mb{Y}_2$ and $\mb{Y}_4$ are not modified, but any transformer side $n$ whose self-admittance is modeled by $\mb{Y}_2$ or $\mb{Y}_4$ contains a constant-impedance load $\mb{Y}_{\mr{L}_n}$ with $\mr{Re}[\mb{Y}_{\mr{L}_n}] \succ \mb{O}$.
\end{corollary}

\subsection{Invertibility using the imaginary parts of $\epsilon'$ and $\epsilon''$}
\label{sec:invertimag}
In simplistic models of  distribution networks, transformer connections are represented by their series reactances, that is, the transformer models assume zero ohmic losses.  In this section, we present an alternative to Theorem \ref{theorem:Yinv} that provides a gateway to guarantee  invertibility of the Y-Bus under such  transformer models. 
Concretely, recall that Theorem \ref{theorem:Yinv} relies on Lemma \ref{lemma:usefulalgebra}. 
% which  proves that matrix $\mb{Y}$ is invertible if   its real part is positive definite and its imaginary part is symmetric, i.e., $\mb{Y}^{-1}$ exists when $\mr{Re}[\mb{Y}] \succ \mb{O}$  $\mr{Im}[\mb{Y}]=\mr{Im}[\mb{Y}]^T$. Notice that 
By applying Lemma \ref{lemma:usefulalgebra} to matrix $j\mb{Y}$,  a dual result can be proved: 
%Matrix $\mb{Y}$ is invertible if its real part is symmetric and its imaginary part is negative definite, i.e., 
$\mb{Y}^{-1}$ exists when $\mr{Re}[\mb{Y}] =\mr{Re}[\mb{Y}]^T$   and $\mr{Im}[\mb{Y}] \prec \mb{O}$ hold. 
The previous result can be used to provide an alternative to Theorem \ref{theorem:Yinv}.  To this end, the following modifications to Assumptions \ref{assumption:open-wye--open-delta} and \ref{assumption:nonidealtransformers} are considered:%\ref{assumption:linecategories}--\ref{assumption:threephases} as well as  \ref{assumption:open-wye--open-delta} and \ref{assumption:missingphases}:
  \renewcommand{\theenumi}{A\arabic{enumi}$''$}
 \renewcommand{\labelenumi}{\theenumi)}
 \begin{enumerate}
\item \label{assumption:insteadtr} Instead of  \ref{assumption:trpd}, it holds that $\mr{Im}[\mb{Z}_{nm}^{-1}] \prec \mb{O}$ and that shunt admittances of the lines are negligible, that is, $\|\mb{Y}_{nm}^{\mr{s}}\| =0$. 
 In \ref{assumption:Y2}--\ref{assumption:Y4} and \ref{assumption:open-wye--open-delta}, the modification of $\mb{Y}_2$ or $\mb{Y}_4$ on self-admittances should be such that $\epsilon'_{i} > 0$ (added  admittance is reactive). \label{assumption:insteadYgDelta}
In  \ref{assumption:Y2} and \ref{assumption:Y4}, the modification of $\mb{Y}_2$ or $\mb{Y}_4$ on mutual admittances should be such that $0 \le \epsilon''_{i} < \epsilon_{i}'$. \label{assumption:insteadyY2}
\setcounter{enumi}{2}
\item Instead of  \ref{assumption:nonidealtransformers}, assume that  $\mr{Im}[y_t] < 0$, that is, transformers have non-zero leakage inductance. \label{assumption:insteadnonidealtr}
 \end{enumerate}

The ensuing Theorem \ref{theorem:Yinvidealtrs} guarantees invertibility of the Y-Bus in distribution circuits where transformers are represented only by their series reactances.  
\begin{theorem}
Under \ref{assumption:insteadtr}, \ref{assumption:connectedness}, \ref{assumption:insteadnonidealtr}, \ref{assumption:missingphases}, $\mb{Y}$ is invertible. \label{theorem:Yinvidealtrs}
\end{theorem}

\begin{IEEEproof}
The proof is similar to that of Theorem \ref{theorem:Yinv}, but use the imaginary operator $\mr{Im}$  instead of $\mr{Re}$,  and replace positive definite arguments with negative definite arguments. 
\end{IEEEproof}

%In the next section, a discussion is provided to guarantee invertibility under more detailed transformer models. 

\subsection{Alternative transformer models}
\label{sec:altxfm}
Nodal admittances provided in Table~\ref{table:transformerConnections}  hold true for a bank of three single-phase transformers.   If no information other than the series admittances for the transformer is available, then the matrices in Table~\ref{table:transformerConnections} provide  a reasonable  model in unbalanced distribution system studies.  This is indeed the case for the feeders we tested. 

 However,  transformers with common-core or shell-type structures are also used in distribution systems.  For such transformer connections,  by  conducting extensive short-circuit measurements, we can first obtain primitive admittance parameters  and then convert them to nodal admittances of the form \eqref{eqngroup:genericSeriesForm} by using  connection matrices \cite{chendillon1974}.  The proof of Section \ref{sec:invertibility} provides a template that shows how each of the block admittances of \eqref{eqngroup:genericSeriesForm} affects the Y-Bus invertibility and can be modified accordingly to the individual transformer model.

As an example, consider the model of a three-legged core-type transformer connected in wye-g--delta given in \cite[eq. (40)]{chendillon1974}. The only difference between that model and the wye-g--delta of Table~\ref{table:transformerConnections} is that the corresponding $\mb{Y}_1$ for the core-type transformer is not diagonal. However, it is mentioned in \cite{chendillon1974} that the off-diagonal elements are considerably smaller in magnitude. Thus, if $\mr{Re}[\mb{Y}_1]$ is diagonally dominant, and hence positive definite,  the wye-g side of the transformer does not compromise the Y-Bus invertibility. Of course,  matrix $\mb{Y}_2$ on the delta side,  still needs to be modified as per \eqref{eqn:Y2prime}. 

 \renewcommand{\theenumi}{\arabic{enumi}}
\section{Numerical results}
\label{sec:numtests}
In this section, we put together the previously developed models to build the Y-Bus matrix and run the Z-Bus power flow for the following distribution  feeders:
\renewcommand{\labelenumi}{(\alph{enumi})}
 \begin{enumerate*}
	\item the IEEE 37-bus feeder,
	\item the IEEE 123-bus feeder,
	\item the  8500-node feeder,
	\item the European 906-bus low voltage feeder (ELV 906-bus feeder).
\end{enumerate*}
The developed MATLAB scripts to conduct the numerical tests  are available at the following page:
\begin{center}
	\texttt{\url{https://github.com/hafezbazrafshan/three-phase-modeling}}
\end{center}
For each feeder, the MATLAB script \texttt{setupYBus} imports the feeder excel files, and builds the required matrices $\mb{Y}_{\mr{net}}$, $\mb{Y}$, and $\mb{Y}_{\mr{NS}}$.  The  script \texttt{setupLoads}  builds the required load parameters for the network according to Table~\ref{table:wyedeltaloads}.   The MATLAB script \texttt{performZBus} performs the iterations of the Z-Bus method given in \eqref{eqn:zbus}. 

The required  data for the feeders are obtained from \cite{ieeefeederdata}.  In order to assess the accuracy of our proposed modeling approach, computed load-flow voltage solutions are compared with benchmark solutions. Specifically, for the IEEE 123-bus and the 8500-node test feeders, the benchmark is obtained from the feeder documents in \cite{ieeefeederdata}. 

The benchmark for the ELV 906-bus feeder  is not provided in \cite{ieeefeederdata}, therefore we use OpenDSS,  the EPRI distribution system simulator \cite{openDSS, openDSSManual}. For the IEEE 37-bus feeder, the benchmark is also obtained from OpenDSS. This is because OpenDSS assumes a relatively high-impedance for the open-delta SVR and thereby allows us to verify the non-ideal SVR models of \eqref{eqngroup:SVRnonIdealYtildes} in Section \ref{sec:svrs}. %\footnote{In order to obtain a benchmark, the original implementation of the IEEE 37-bus and the ELV 906-bus feeder  as provided by the OpenDSS software were not modified.}

The original  data  from \cite{ieeefeederdata} do not provide the series impedance for any of the SVRs in any of the feeders.   Unlike the OpenDSS model of the IEEE 37-bus feeder,    OpenDSS models of IEEE 123-bus and the 8500-node feeders  assume  very small series impedances for  SVRs.  For the IEEE 123-bus and the 8500-node feeders, it will be demonstrated that using the ideal SVR models in \eqref{eqngroup:SVRIdealYtildes} does not affect the accuracy of power flow solutions. 

\begin{figure}[t]
	\centering 
\subfloat[]{\includegraphics[scale=0.22]{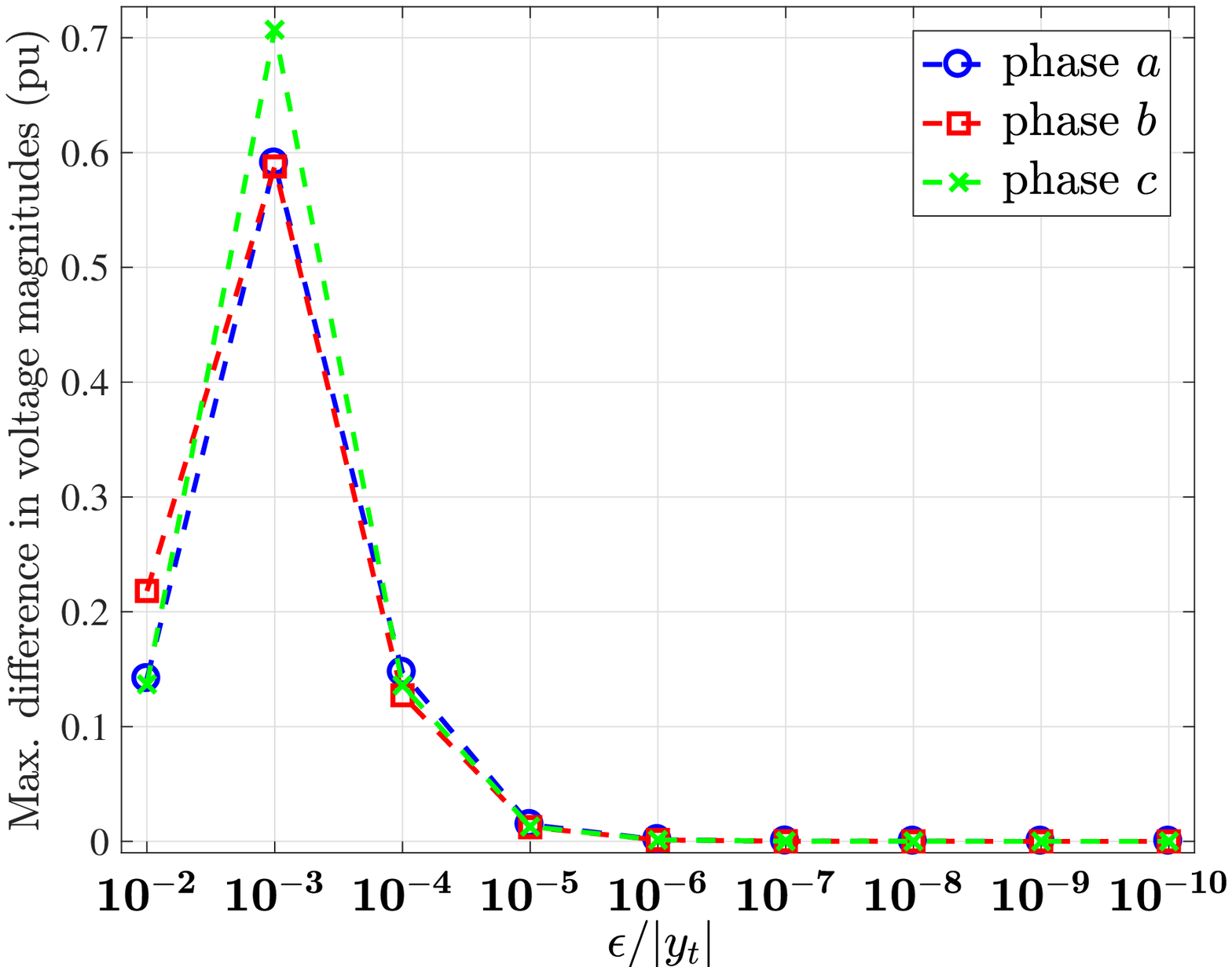}\label{fig:IEEE37epsilon} } 
\subfloat[]{\includegraphics[scale=0.22]{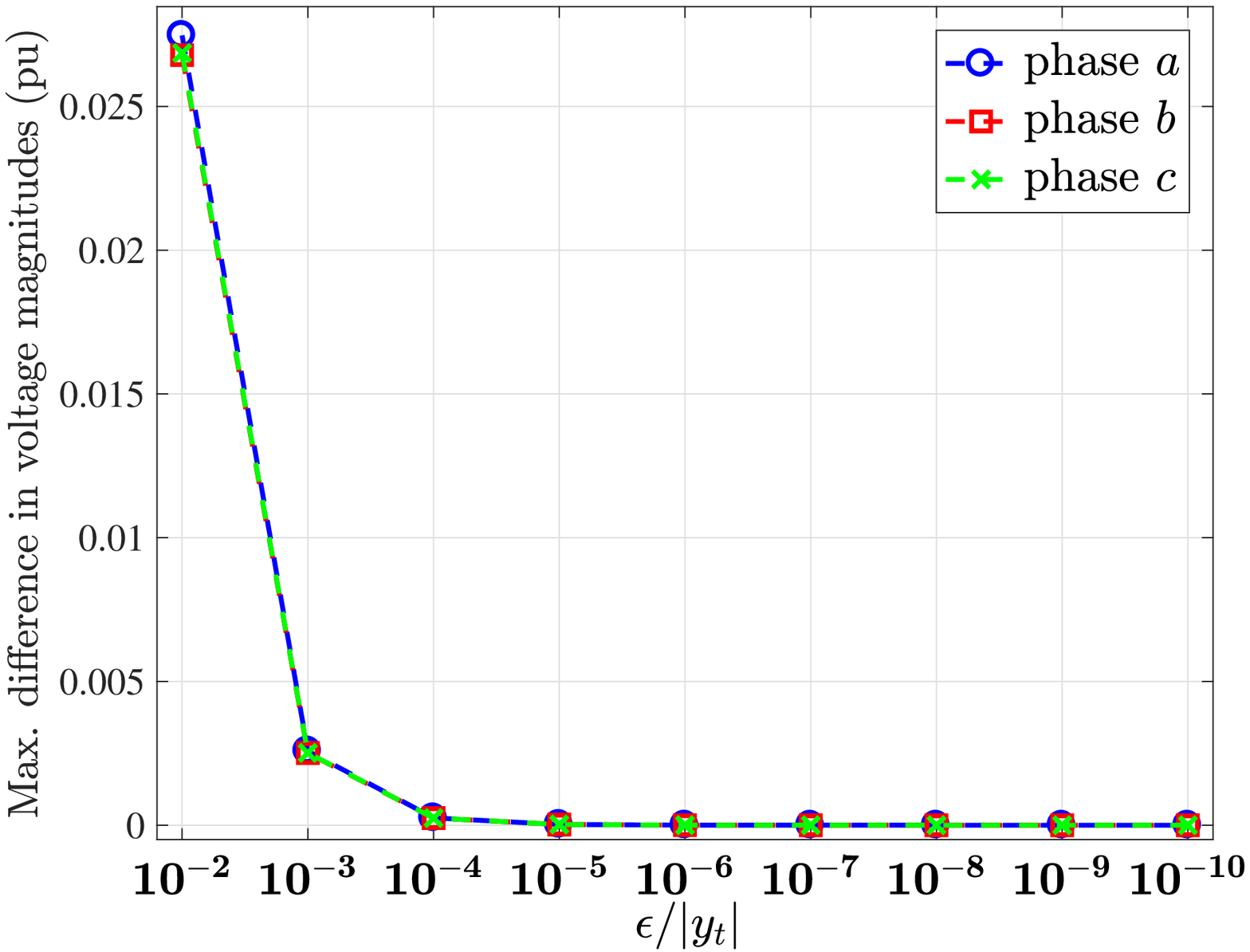}  \label{fig:IEEE123epsilon}}
	\caption{Maximum difference in computed voltage magnitudes obtained from  load-flow runs on \protect\subref{fig:IEEE37epsilon} the IEEE 37-bus feeder and \protect\subref{fig:IEEE123epsilon} the IEEE 123-bus feeder. The computed voltage solutions in successive load-flow runs converge when the value of $\epsilon/|y_t|$ reaches below $10^{-5}$ and $10^{-3}$ for the IEEE 37-bus feeder and the IEEE 123-bus feeder respectively. 	\label{fig:IEEE37123epsilon} }
\end{figure}

In the modeling of the delta--delta transformers in IEEE 37-bus and IEEE 123-bus feeders, we use nodal admittances from Table~\ref{table:sixcategories}. By setting $\epsilon'=2\epsilon''=\epsilon$,   we perform a series of Z-Bus power flow runs with values of $\epsilon$ ranging from $10^{-2}|y_t|$ to $10^{-10}|y_t|$.  The maximum difference in voltage magnitudes for successive values of $\epsilon$ is depicted in Fig.~\ref{fig:IEEE37123epsilon}\protect\subref{fig:IEEE37epsilon} for the IEEE 37-bus feeder and in Fig.~\ref{fig:IEEE37123epsilon}\protect\subref{fig:IEEE123epsilon}  for the IEEE 123-bus feeder. It turns out that when $\epsilon$ becomes smaller than a threshold, the computed voltages from successive power flow runs converge. This  threshold as a fraction of $\epsilon / |y_t|$ is  respectively $10^{-5}$ and $10^{-3}$ for the IEEE 37-bus and the IEEE 123-bus feeders. The IEEE 37-bus contains two delta--delta transformers whereas IEEE 123-bus includes only one, therefore the required threshold for $\epsilon$ is smaller in the IEEE 37-bus than that of IEEE 123-bus.

\begin{figure}[t]
 	\centering 
 	\subfloat[]{\includegraphics[scale=0.22]{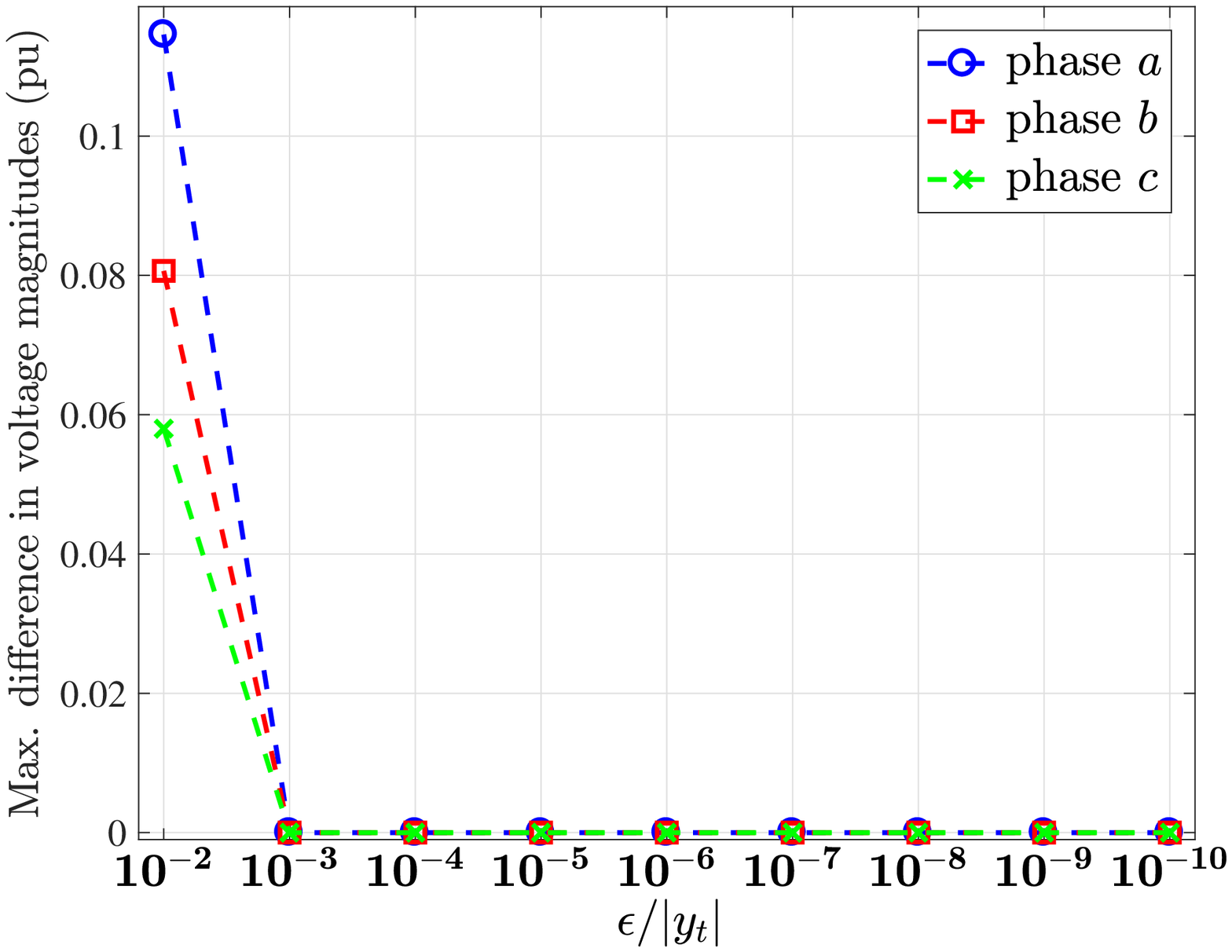}\label{fig:IEEE8500epsilon} } 
 	\subfloat[]{ \includegraphics[scale=0.22]{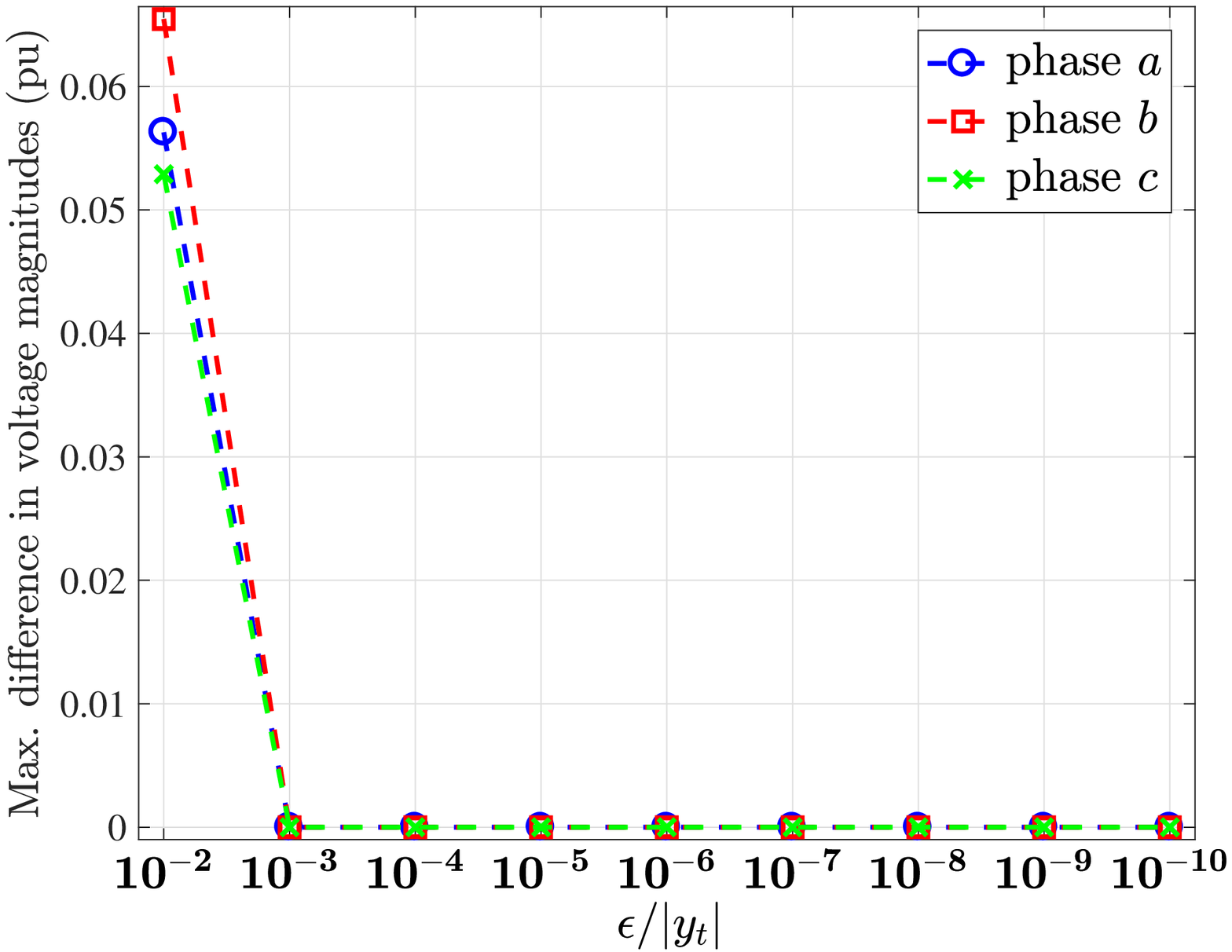}  \label{fig:IEEELVepsilon}}
 	\caption{Maximum difference in computed voltage magnitudes obtained from  load-flow runs on \protect\subref{fig:IEEE8500epsilon} the  8500-node feeder and \protect\subref{fig:IEEELVepsilon} the ELV 906-bus feeder. The computed voltage solutions in successive load-flow runs converge when the value of $\epsilon/|y_t|$ reaches below  $10^{-3}$ for both feeders. 	\label{fig:IEEE8500LVepsilon} }
 \end{figure}

In the modeling of the delta--wye transformers in  8500-node and the ELV 906-bus feeders,  we use nodal admittances of Table~\ref{table:sixcategories} with $\epsilon'=\epsilon$. Once again, a series of load-flow computations are performed with values of $\epsilon$ ranging from $10^{-2}|y_t|$ to $10^{-10}|y_t|$.  The maximum difference in voltage magnitudes for various values of $\epsilon$ are depicted in Fig.~\ref{fig:IEEE8500LVepsilon}\protect\subref{fig:IEEE8500epsilon} for the 8500-node feeder and in Fig.~\ref{fig:IEEE8500LVepsilon}\protect\subref{fig:IEEELVepsilon} for the ELV 906-bus feeder. Once the value of  $\epsilon / |y_t|$ reaches below $10^{-3}$, the voltage solutions converge.  In all four feeders, if we do not modify $\mb{Y}_2$, the resultant Y-Bus is singular. 
 
Table \ref{table:summary} summarizes the comparison of our results with  benchmark solutions using the fixed value of $\epsilon=10^{-6}|y_t|$. The reported values are below $0.75\%$ of the benchmark solutions. Description of each test feeder is provided next, along with calculated voltage plots corresponding to Table \ref{table:summary}.  In the voltage plots,  bus labels of each feeder have been modified to represent successive numbers starting from 1.  
 \begin{table}[t]
 	\centering
 	\caption{Maximum $\mr{pu}$ difference in computed voltage magnitudes vs. benchmark}
 	\begin{tabular}{c|c|c|c}
 		Feeder  & phase $a$ & phase $b$ & phase $c$  \\
 		\hline
 		\hline
 		IEEE 37-bus  & 0.0067  &  0.0019  &  0.0053 \\
 		\hline
 		IEEE 123-bus  &0.0061&0.0034&0.0039  \\
 		\hline 
 		8500-node  &0.0010&0.0013&0.0034\\
 		\hline 
 		ELV 906-bus  & 0.0048&0.0055&0.0021 \\
 		\hline
 	\end{tabular}
 	\label{table:summary}
 \end{table}

\subsection{IEEE 37-bus feeder}
\begin{figure}[t]
	\centering
	\includegraphics[scale=0.50]{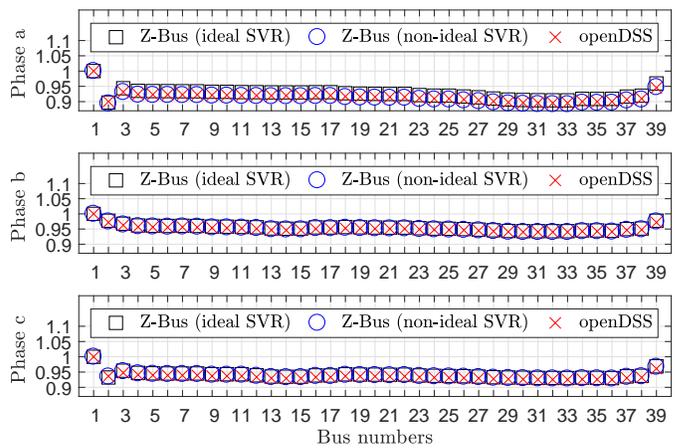}
	\caption{The IEEE 37-bus  feeder voltage profile obtained from the Z-Bus method using ideal SVRs (black squares) and non-ideal SVRs (blue circles) in comparison to the voltage profile provided by OpenDSS (red cross). Bus labels have been modified to represent successive numbers starting from 1.} \label{fig:IEEE37magnitudes}
\end{figure}

%\begin{figure}[t]
%	\centering
%	\subfloat[]{\includegraphics[scale=0.50]{Figures/voltage_plots/IEEE37magnitudes.eps} \label{fig:IEEE37magnitudes}} 
%	\hfill
%	\subfloat[]{\includegraphics[scale=0.50]{Figures/voltage_plots/IEEE37phases.eps} \label{fig:IEEE37phases}}
%	\caption{ IEEE-37 voltage profile obtained from the Z-Bus method using ideal regulators (black squares) and non-ideal regulators (blue circles) in comparison to the IEEE-37 voltage profile provided by OpenDSS (red cross):  \protect\subref{fig:IEEE37magnitudes} voltage magnitudes \protect\subref{fig:IEEE37phases} voltage phases. Bus labels have been modified to represent successive numbers starting from 1.}
%	\label{fig:IEEE37voltageprofile}
%\end{figure}

\begin{table}[t]
	\centering
	\caption{Ideal SVRs vs Non-ideal SVRs: Max. $\mr{pu}$ difference in computed voltage magnitudes vs. OpenDSS}
	\begin{tabular}{c|c|c|c}
		SVR model & phase $a$ & phase $b$ & phase $c$  \\
		\hline
		\hline
		Ideal [eq. \eqref{eqngroup:SVRIdealYtildes}] &0.0136&0.0046&0.0072 \\
		\hline
		Non-ideal [eq. \eqref{eqngroup:SVRnonIdealYtildes}] & 0.0067  &  0.0019  &  0.0053 \\
		\hline
	\end{tabular}
	\label{table:idealvsnonidealIEEE37}
\end{table}
\label{sec:numtests:IEEE37}
The IEEE 37-bus feeder features the following rather distinctive characteristics: \renewcommand{\labelenumi}{(\alph{enumi})}
\begin{enumerate*}
	\item a delta--delta substation transformer, rated $2500~\mr{kVA}$,  $230~\mr{kV}/4.8~\mr{kV}$ line to line, with  $z_t=(2+j8)\%$,  on edge $(1,2)$;
	\item a delta--delta transformer, rated $500~\mr{kVA}$,  $4.8~\mr{kV}/0.48~\mr{kV}$ line to line, with $z_t=0.09+j1.81 \%$, on edge $(24,38)$;
	\item an open-delta SVR   on edge $(2,3)$, with a relatively high impedance of $z_t=j1 \%$; and 
	\item a variety of delta-connected constant-power, constant-current, and constant-impedance loads.
\end{enumerate*}

 The impedance of the SVR was obtained from the OpenDSS implementation provided in the file \texttt{\url{OpenDSS/Distrib/IEEETestCases/37Bus/ieee37.dss}}.   Voltage magnitudes corresponding to the value  of $\epsilon=10^{-6}|y_t|$ are plotted in Fig.~\ref{fig:IEEE37magnitudes} for both ideal and non-ideal regulator models where bus $39$ represents the $n'$ of the SVR.  The solution provided by OpenDSS is also provided for verification.  Phase plots have been omitted here due to space limitations but are available on our github page. The maximum voltage magnitude  difference  between the computed solutions and those of OpenDSS are tabulated in Table~\ref{table:idealvsnonidealIEEE37}.  It is inferred from Table \ref{table:idealvsnonidealIEEE37} that when the series impedance of the SVR is relatively high, the non-ideal SVR models in \eqref{eqngroup:SVRnonIdealYtildes} provide more accurate results.

%\begin{figure*}[t]
%	\centering
%	\subfloat[]{\includegraphics[scale=0.40]{Figures/voltage_plots/IEEE123magnitudes.eps} \label{fig:IEEE123magnitudes}} 
%	\hfill
%	%	\subfloat[]{\includegraphics[scale=0.40]{Figures/voltage_plots/IEEE123phases.eps} \label{fig:IEEE123phases}}
%	\caption{ IEEE-123 voltage profile obtained from the Z-Bus method using ideal regulators (black squares) and non-ideal regulators (blue circles) in comparison to the IEEE-123 voltage profile provided by the benchmark (red cross):  \protect\subref{fig:IEEE123magnitudes} voltage magnitudes \protect\subref{fig:IEEE123phases} voltage angles. Bus labels have been modified to represent successive numbers starting from 1. }
%	\label{fig:IEEE123voltageprofile}
%\end{figure*}

\subsection{IEEE 123-bus test feeder}
\begin{figure}[t]
	\centering
	\includegraphics[scale=0.50]{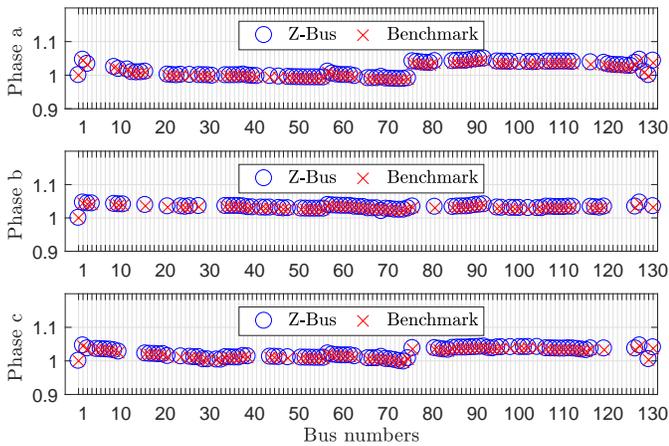} 
	\hfill
	\caption{The IEEE 123-bus feeder voltage profile obtained from the Z-Bus method  (blue circles) in comparison to the voltage profile provided by the benchmark (red cross). Bus labels have been modified to represent successive numbers starting from 1. }
	\label{fig:IEEE123magnitudes}
\end{figure}

The IEEE 123-bus test feeder features
\begin{enumerate*}
	\item three-, two-, and single-phase laterals;
	\item four wye-connected SVRs, namely, ID1 on edge $(1,2)$ three-phase gang-operated, ID2 on phase $a$ of edge $(12,13)$, ID3 on phases $a,c$ of edge $(28,29)$, ID4 on  three-phases of edge $(75,76)$;
	\item a delta-delta transformer, rated $150~\mr{kVA}$ and $4.16~\mr{kV}/0.48 ~\mr{kV}$ line to line, on edge $(68,69)$.
	\end{enumerate*}
The voltage profile is provided in Fig.~\ref{fig:IEEE123magnitudes}, where the buses $127, 128, 129,130$ are the $n'$ of SVRs with IDs $1,2,3,4$ respectively (the last four markers on the plot). The corresponding phase plots  are available on our github page. 

\subsection{8500-node test feeder}
The 8500-node test feeder includes approximately $2500$ medium voltage buses\footnote{The feeder documentation distinguishes between the term ``bus" and ``node."  A bus can have multiple nodes (phases). The feeder includes $2500$ primary buses, each of which may be three-, two-, or one-phase.}  \cite{IEEE8500}.   The network features
\begin{enumerate*}
	\item three-, two-, and single-phase laterals;
	\item a delta--wye substation transformer on edge $(1,2)$, rated $27.5~\mr{MVA}$,  $115~\mr{kV}/12.47~\mr{kV}$ line to line, with  $z_t=1.344+j15.51 \%$; and
		\item four wye-connected, individually operated, three-phase SVRs, namely, ID1 on edge $(2,3)$, ID2 on edge $(201,202)$, ID3 on edge $(146,147)$, ID4 on edge $(1777,1778)$.
\end{enumerate*}

\begin{figure*}[t]
	\centering
	\subfloat[]{\includegraphics[scale=0.40]{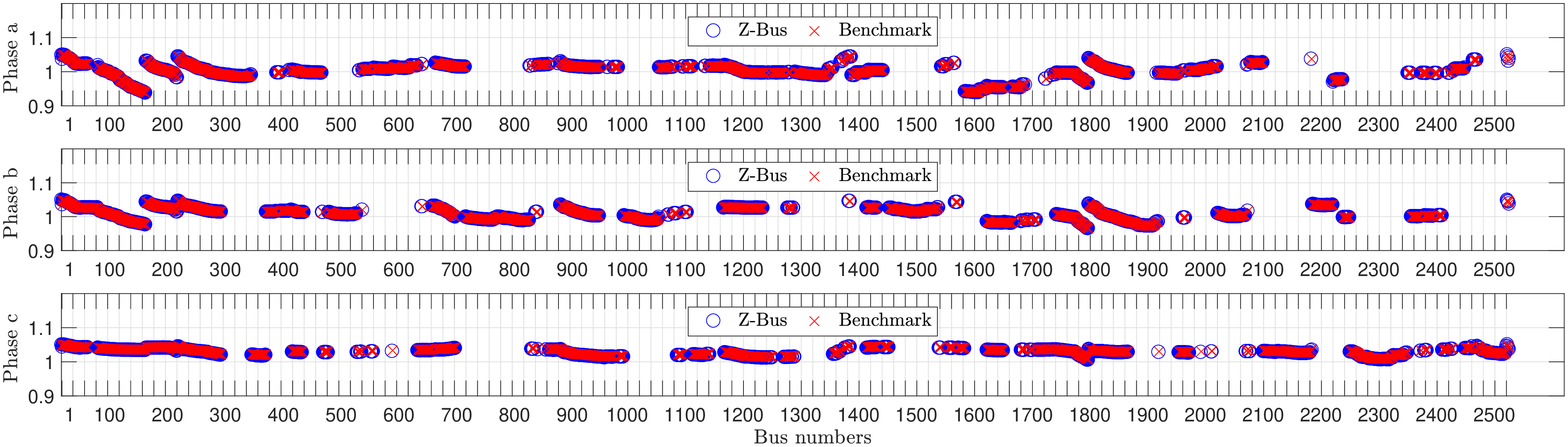} \label{fig:IEEE8500magnitudes}} 
	\hfill
	\subfloat[]{\includegraphics[scale=0.40]{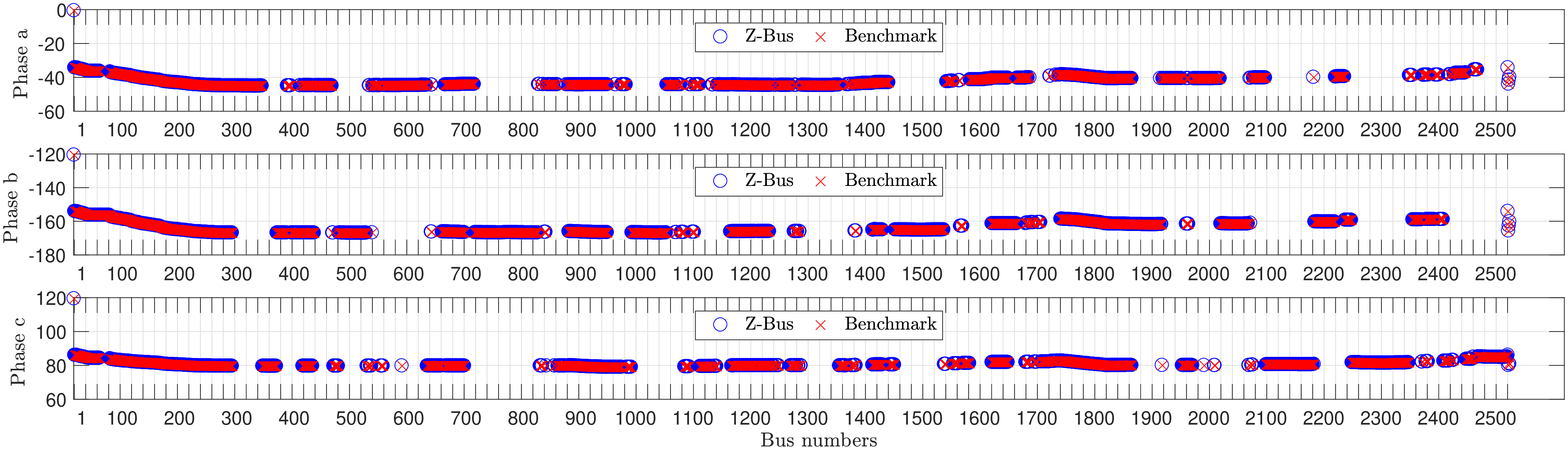} \label{fig:IEEE8500phases}}
	\caption{The  8500-node feeder voltage profile obtained from the Z-Bus method using ideal regulators (blue circles) in comparison to the  voltage profile provided by the benchmark (red cross): \protect\subref{fig:IEEE8500magnitudes} voltage magnitudes \protect\subref{fig:IEEE8500phases} voltage phases.}
	\label{fig:IEEE8500voltageprofile}
\end{figure*}

The voltage profile is provided in Fig.~\ref{fig:IEEE8500voltageprofile}, where the buses $2502, 2503, 2504, 2505$ are the $n'$ of SVRs with IDs $1,2,3,4$ respectively (the last four markers on the plot). For this feeder, the SVR models of \eqref{eqngroup:SVRIdealYtildes} should correspond to type A SVRs. In this case, the $\mp$ sign in \eqref{eqn:1phaseregtap} is changed to $\pm$, and instead of $\mb{A}_{\mr{i}}$, its inverse $\mb{A}_{\mr{i}}^{-1}$ is used. 

\subsection{The European 906-bus low voltage feeder}
\begin{figure*}[t]
	\centering
	\includegraphics[scale=0.40]{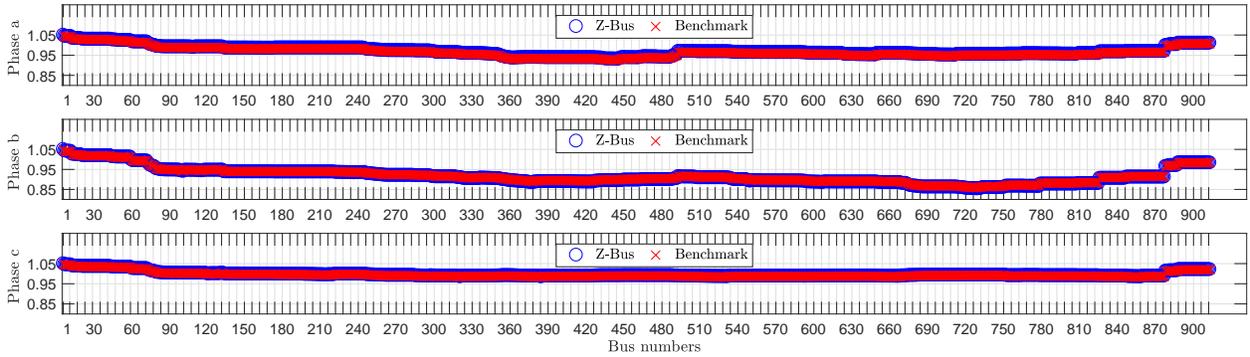} 
	\hfill
	\caption{ The ELV 906-bus feeder voltage profile obtained from the Z-Bus method  (blue circles) in comparison to the  voltage profile provided by the benchmark (red cross). The intensity of the plot is due to the fact that this feeder contains three-phase laterals only. \label{fig:IEEELVmagnitudes}}
\end{figure*}

The ELV feeder features \begin{enumerate*}
	\item three-phase laterals (there are no missing phases), and  
	\item a delta--wye substation transformer, rated $800~\mr{kVA}$, $11~\mr{kV}/416~\mr{V}$ line to line, with $z_t=(0.4+j4)\%$, on edge $(1,2)$. 
\end{enumerate*}
The obtained voltage profile is provided in Fig.~\ref{fig:IEEELVmagnitudes}. The  phase plots are omitted due to space limitation, but are available on our github page. 

\section{Concluding remarks}
\label{sec:conclusion}
This paper focused on nodal admittance modeling of three-phase distribution networks. Models for transmission lines and most relevant transformer connections are reviewed, while novel models for step-voltage regulators are derived, explicitly accounting for their tap positions and their series admittance. Putting the models for the series elements together yields the network bus admittance matrix Y-Bus, whose invertibility is crucial in deriving power flow solutions via the Z-Bus method. The paper carefully laid out the conditions on each series element that guarantee the invertibility of the Y-Bus and proves why previous proposals on modifications of certain transformer connections restore its invertibility. The conditions are tailored to practical distribution networks, which can be radial or meshed, feature any number of transformers and SVRs, and include missing phases. Comprehensive numerical tests are presented for the IEEE 37-bus, the IEEE 123-bus, the 8500-node medium-voltage feeders, and the European 906-bus low-voltage feeder. The codes that build the Y-Bus and compute the power flow solutions are provided online.

\appendices
\section{Useful linear algebra results}
\label{sec:appendixLemma1}
\begin{lemma}
\label{lemma:usefulalgebra}
Consider matrix $\mb{Y}=\mb{G}+ j\mb{B} \in \mbb{C}^{n \times n}$  with  $\mb{G} \succ \mb{O}$ and $\mb{B}=\mb{B}^T$. Then $\mb{Y}^{-1}$ exists and $\mr{Re}[{\mb{Y}^{-1}}]\succ \mb{O}$. 
\end{lemma} 
\begin{IEEEproof} The invertibility proof is based on contradiction.  Assume $\mb{x}= \mb{x}_{\mr{R}} + j \mb{x}_{\mr{I}} \in \mbb{C}^{n}$ is a nonzero vector in the nullspace of $\mb{Y}$. Then, we have that  
		\begin{IEEEeqnarray}{rCl}
			\mb{Y}\mb{x} &=& (\mb{G}+ j \mb{B}) (\mb{x}_{\mr{R}} + j \mb{x}_{\mr{I}} ) \notag \\ 
			&=& \mb{G} \mb{x}_{\mr{R}}- \mb{B} \mb{x}_{\mr{I}} + j (\mb{G}\mb{x}_{\mr{I}} + \mb{B}\mb{x}_{\mr{R}}) = \mb{0} \label{eqn:Yxzero}
		\end{IEEEeqnarray}
Thus, we conclude that 
\begin{IEEEeqnarray}{rClrCl}
			\mb{G}\mb{x}_{\mr{R}} - \mb{B}\mb{x}_{\mr{I}} = \mb{0},  \quad& 
			\mb{G}\mb{x}_{\mr{I}} + \mb{B} \mb{x}_{\mr{R}} = \mb{0}. \label{eqn:gxi-bxr}
\end{IEEEeqnarray}
Multiplying the second equality in \eqref{eqn:gxi-bxr} by $\mb{x}_{\mr{I}}^T$ yields 
$\mb{x}_{\mr{I}}^T \mb{G} \mb{x}_{\mr{I}} + \mb{x}_{\mr{I}}^T \mb{B} \mb{x}_{\mr{R}} = 0.  $
In the latter equality, the first term is non-negative due to the assumption that  $\mb{G} \succ \mb{O}$ . Therefore we must have that 
	$\mb{x}_{\mr{I}}^T \mb{B} \mb{x}_{\mr{R}} \le 0$.  Since $\mb{B}= \mb{B}^T$ we find that 
	\begin{IEEEeqnarray}{rCl}
		\mb{x}_{\mr{I}}^T \mb{B} \mb{x}_{\mr{R}} = \mb{x}_{\mr{R}}^T \mb{B}^T \mb{x}_{\mr{I}} =  \mb{x}_{\mr{R}}^T \mb{B}\mb{x}_{\mr{I}} \le 0. \label{eqn:xrtbxi}
	\end{IEEEeqnarray}
	
Multiplying the first equality in \eqref{eqn:gxi-bxr} by $\mb{x}_{\mr{R}}^T$ gives $\mb{x}_{\mr{R}}^T \mb{G} \mb{x}_{\mr{R}} - \mb{x}_{\mr{R}}^T \mb{B} \mb{x}_{\mr{I}} = 0.$  In the latter equality, the first term is non-negative due to the assumption that $\mb{G} \succ \mb{O}$, and the second term is non-negative due to \eqref{eqn:xrtbxi}.  Therefore, both terms must be zero.  Hence $\mb{x}_{\mr{R}}= \mb{0}$. Replacing $\mb{x}_{\mr{R}} = \mb{0}$ in the second equality of \eqref{eqn:gxi-bxr} yields $\mb{x}_{\mr{I}}= \mb{0}$ as well. Hence, $\mb{x}= \mb{0}$ which is a contradiction, since  the vector $\mb{x}$ was assumed to be nonzero. 	
To prove the second portion, let 	$\mb{Y}^{-1}=\mb{R}+j\mb{X}$. Then, it holds that $(\mb{R}+j\mb{X}) (\mb{G}+j\mb{B}) = \mb{I}$. Hence, 
		\begin{subequations}
			\label{eqngroup:RGXB}
			\eq{rCl}{
				\mb{R}\mb{G}-\mb{X}\mb{B} = \mb{I} \label{eqn:RGXB} \\
				\mb{X}\mb{G}+ \mb{R}\mb{B}= \mb{O}. \label{eqn:XGRB}
			}
			Using \eqref{eqn:XGRB}, we obtain $\mb{X}= -\mb{R} \mb{B} \mb{G}^{-1}$. The latter expression can be inserted into \eqref{eqn:RGXB} to obtain:
			\eq{rCl}{
				\mb{R}\mb{G}+ \mb{R} \mb{B} \mb{G}^{-1} \mb{B} = \mb{I} \Rightarrow  \mb{R} = (\mb{G}+ \mb{B} \mb{G}^{-1} \mb{B})^{-1}.  \label{eqn:Rinv}
			}
			Notice that $\mb{G}+ \mb{B} \mb{G}^{-1} \mb{B}$  is positive definite since $\mb{G} \succ \mb{O}$. Hence, it holds that $\mb{R} \succ \mb{O}$. 
		\end{subequations}
\end{IEEEproof}

\color{black}
\begin{lemma} \label{lemma:utWv}
	For vectors $\mb{u}$ and $\mb{v}$ and matrix $\mb{W}$ we have that
	\small{\eq{rCl}{\mb{u}^T \mb{W}  \mb{u}  - \mb{u}^T\mb{W}  \mb{v}  +  \mb{v}^T \mb{W}\mb{v} - \mb{v}^T \mb{W}\mb{u}
		\hspace{-0.05cm}=  \hspace{-0.05cm} (\mb{u}-\mb{v})^T \mb{W} (\mb{u}-\mb{v}).\IEEEeqnarraynumspace \label{eqn:utWu} }}
\end{lemma}
\begin{IEEEproof}
	The left-hand side of \eqref{eqn:utWu} can be factored
	\small{\eq{+rCl+x*}{ && \mb{u}^T (\mb{W} \mb{u}- \mb{W} \mb{v})  - \mb{v}^T ( \mb{W} \mb{u}- \mb{W} \mb{v})= (\mb{u}-\mb{v})^T \mb{W} (\mb{u}-\mb{v}). \nonumber  & \IEEEQEDhere }}
\end{IEEEproof}

\begin{lemma} \label{lemma:utWXZv}
	For vectors $\mb{u}$, $\mb{v}$ and matrices $\mb{W}$, $\mb{X}$,  $\mb{Z}$  we have 
\eq{rCl}{\mb{u}^T \mb{W}  \mb{u}  - \mb{u}^T\mb{X}  \mb{v}  +  \mb{v}^T \mb{Z}\mb{v} - \mb{v}^T \mb{X}^T\mb{u} = \notag  \\ 
		\mspace{100mu}  \bmat{ \mb{u}^T & \mb{v}^T} \bmat{ \mb{W} & - \mb{X} \\ - \mb{X}^T & \mb{Z}} \bmat{\mb{u} \\\mb{v} } \label{eqn:utWXZv}}
\end{lemma}
\begin{IEEEproof}
	The right-hand side of \eqref{eqn:utWXZv} is written as follows
	\small{\eq{+rCl+x*}{&&\bmat{\mb{u} \\ \mb{v}}^T \bmat{ \mb{W} \mb{u} - \mb{X} \mb{v} \\ - \mb{X}^T \mb{u} + \mb{Z} \mb{v}} 
		=\mb{u}^T \mb{W}  \mb{u}  - \mb{u}^T\mb{X}  \mb{v} - \mb{v}^T \mb{X}^T\mb{u} +  \mb{v}^T \mb{Z}\mb{v}. \nonumber & \IEEEQEDhere }}
\end{IEEEproof}

%\begin{lemma} \label{lemma:Ynmninvertible}
%	For a  transmission line $(n,m)$,  $\mb{Y}_{nm}^{(n)}$ is symmetric, invertible, and $\mr{Re}[(\mb{Y}_{nm}^{(n)})^{-1}]$ is positive definite. 
%	\begin{IEEEproof}
%	Due to \eqref{eqngroup:txlineactual} and assumption \ref{assumption:trsymmetric}, $\mb{Y}_{nm}^{(n)}$ is clearly symmetric.  Based  on 
% assumption \ref{assumption:trpd} and Lemma \ref{lemma:usefulalgebra}, $\mb{Y}_{nm}^{(n)}$ is invertible and $\mr{Re}[(\mb{Y}_{nm}^{(n)})^{-1}] \succ \mb{O}$. 
%	\end{IEEEproof}
%\end{lemma}

\begin{lemma} \label{lemma:Fr} Suppose Assumption~\ref{assumption:transmissionlines} of Section~\ref{sec:invertibility} holds for the line $(n',m)$ of the SVR. For the values of $\mb{A}_{\mr{i}}$ and $\mb{Z}_{\mr{R}}$ in Table~\ref{table:vrmatrices},  the matrix $\mb{F}_{\mr{R}} = \mb{I}_{|\Omega_n|} + \mb{Y}_{n'm}^{(n')} \mb{A}_{\mr{i}}^T\mb{Z}_{\mr{R}} \mb{A}_{\mr{i}}$ [cf.~\eqref{eqn:Freg}] is invertible. Moreover, the following holds for $\mb{F}_{\mr{R}}$:
\eq{rCl}{\mb{F}_{\mr{R}}^{-T}= \mb{I}_{|\Omega_n|} -\mb{A}_{\mr{i}}^T \mb{Z}_{\mr{R}} \mb{A}_{\mr{i}} \mb{F}_{\mr{R}}^{-1} \mb{Y}_{n'm}^{(n')}. \label{eqn:Frinv}}
\end{lemma}
\begin{IEEEproof}
Combining~\eqref{eqngroup:transmissionLineYtildes},  Assumption~\ref{assumption:transmissionlines} of Section~\ref{sec:invertibility}, and Lemma~\ref{lemma:usefulalgebra}, it follows that $\mb{Y}_{nm}^{(n)}$ is symmetric and invertible, and $\mr{Re}[(\mb{Y}_{nm}^{(n)})^{-1}] \succ \mb{O}$.  Therefore, we can write 
\eq{rCl}{
\mb{F}_{\mr{R}} = \mb{Y}_{n'm}^{(n')} \left[  (\mb{Y}_{n'm}^{(n')})^{-1} + \mb{A}_{\mr{i}}^T \mb{Z}_{\mr{R}} \mb{A}_{\mr{i}} \right]. \label{eqn:FSVRproduct}
}
Hence, $\mb{F}_{\mr{R}}$ is invertible if and only if $ (\mb{Y}_{n'm}^{(n')})^{-1} + \mb{A}_{\mr{i}}^T \mb{Z}_{\mr{R}} \mb{A}_{\mr{i}}$ is invertible.  
For all SVRs, $\mb{Z}_{\mr{R}}$ is diagonal and $\mr{Re}[\mb{Z}_{\mr{R}}]\succeq \mb{O}$ holds due to ohmic or zero impedances.  Furthermore, the gain matrix  $\mb{A}_{\mr{i}}$ is real-valued and  $ \mb{A}_{i}^T \mb{Z}_{\mr{R}} \mb{A}_{\mr{i}}$ is symmetric. Also recall that $\mr{Re}[(\mb{Y}_{nm}^{(n)})^{-1}] \succ \mb{O}$.   Therefore, it holds that $ \mr{Re}[(\mb{Y}_{n'm}^{(n')})^{-1} + \mb{A}_{\mr{i}}^T \mb{Z}_{\mr{R}} \mb{A}_{\mr{i}}] \succ \mb{O}$.  It follows from Lemma~\ref{lemma:usefulalgebra} that $ (\mb{Y}_{n'm}^{(n')})^{-1} + \mb{A}_{\mr{i}}^T \mb{Z}_{\mr{R}} \mb{A}_{\mr{i}}$ is invertible. This proves $\mb{F}_{\mr{R}}^{-1}$ exists.

%Next, take the inverse of \eqref{eqn:FSVRproduct} and multiply both sides by $\mb{Y}_{n'm}^{(n')}$: 
It follows from~\eqref{eqn:FSVRproduct} that
\eq{rCl}{
 \mb{F}_{\mr{R}}^{-1}\mb{Y}_{n'm}^{(n')} &=& \left[(\mb{Y}_{n'm}^{(n')})^{-1} + \mb{A}_{\mr{i}}^T \mb{Z}_{\mr{R}} \mb{A}_{\mr{i}}\right]^{-1} \label{eqn:FrYnmInv}.
}
Due to the symmetricity of $(\mb{Y}_{n'm}^{(n')})^{-1}$ and $\mb{A}_{\mr{i}}^T \mb{Z}_{\mr{R}} \mb{A}_{\mr{i}}$, it follows from \eqref{eqn:FrYnmInv} that $\mb{F}_{\mr{R}}^{-1}\mb{Y}_{n'm}^{(n')}$ is symmetric:
\eq{rCl}{\mb{F}_{\mr{R}}^{-1}\mb{Y}_{n'm}^{(n')} &=& (\mb{Y}_{n'm}^{(n')})\mb{F}_{\mr{R}}^{-T}. \label{eqn:FrTranspose}}

Multiplying  \eqref{eqn:Freg}  from the left by  $\mb{F}_{\mr{R}}^{-1}$ yields: 
\eq{rCl}{
\mb{I}_{|\Omega_n|} &=& \mb{F}_{\mr{R}}^{-1} +\mb{F}_{\mr{R}}^{-1} \mb{Y}_{n'm}^{(n')} \mb{A}_{\mr{i}}^T\mb{Z}_{\mr{R}} \mb{A}_{\mr{i}}  \label{eqn:Frinv1} \IEEEeqnarraynumspace \\
%\Rightarrow   \mb{F}_{\mr{R}}^{-1} &=& \mb{I}_{|\Omega_n|}  - \mb{F}_{\mr{R}}^{-1} \mb{Y}_{n'm}^{(n')} \mb{A}_{\mr{i}}^T\mb{Z}_{\mr{R}} \mb{A}_{\mr{i}}  \label{eqn:Frinv2} \\
\Rightarrow  \mb{F}_{\mr{R}}^{-1} &=& \mb{I}_{|\Omega_n|}  -  \mb{Y}_{n'm}^{(n')} \mb{F}_{\mr{R}}^{-T} \mb{A}_{\mr{i}}^T\mb{Z}_{\mr{R}} \mb{A}_{\mr{i}} 
 \label{eqn:Frinv3},}
where in the last equality we replaced $ \mb{F}_{\mr{R}}^{-1} \mb{Y}_{n'm}^{(n')}$ with its equivalent from \eqref{eqn:FrTranspose}. 
Transposing \eqref{eqn:Frinv3} yields \eqref{eqn:Frinv}. 
\end{IEEEproof}
 \color{black}

\ifCLASSOPTIONcaptionsoff
  \newpage
\fi

% bibliography
\bibliographystyle{IEEEtran}
{
\bibliography{bibliography}

\begin{thebibliography}{10}
\providecommand{\url}[1]{#1}
\csname url@samestyle\endcsname
\providecommand{\newblock}{\relax}
\providecommand{\bibinfo}[2]{#2}
\providecommand{\BIBentrySTDinterwordspacing}{\spaceskip=0pt\relax}
\providecommand{\BIBentryALTinterwordstretchfactor}{4}
\providecommand{\BIBentryALTinterwordspacing}{\spaceskip=\fontdimen2\font plus
\BIBentryALTinterwordstretchfactor\fontdimen3\font minus
  \fontdimen4\font\relax}
\providecommand{\BIBforeignlanguage}[2]{{%
\expandafter\ifx\csname l@#1\endcsname\relax
\typeout{** WARNING: IEEEtran.bst: No hyphenation pattern has been}%
\typeout{** loaded for the language `#1'. Using the pattern for}%
\typeout{** the default language instead.}%
\else
\language=\csname l@#1\endcsname
\fi
#2}}
\providecommand{\BIBdecl}{\relax}
\BIBdecl

\bibitem{BirtGraffyMcDonalElAbiad1976}
K.~A. Birt, J.~J. Graffy, J.~D. McDonald, and A.~H. El-Abiad, ``Three phase
  load flow program,'' \emph{IEEE Trans. Power App. Syst.}, vol.~95, no.~1, pp.
  59--65, Jan. 1976.

\bibitem{GarciaPereiraCarneiroCostaMartins2000}
P.~A.~N. Garcia, J.~L.~R. Pereira, S.~Carneiro, V.~M. da~Costa, and N.~Martins,
  ``Three-phase power flow calculations using the current injection method,''
  \emph{IEEE Trans. Power Syst.}, vol.~15, no.~2, pp. 508--514, May 2000.

\bibitem{Chen1991pf}
T.~H. Chen, M.~S. Chen, K.~Hwang, P.~Kotas, and E.~A. Chebli, ``Distribution
  system power flow analysis - a rigid approach,'' \emph{IEEE Trans. Power
  Del.}, vol.~6, no.~3, pp. 1146--1152, July 1991.

\bibitem{AraujoPenidoCarneiroPereira2013}
L.~R. Araujo, D.~R.~R. Penido, S.~Carneiro, and J.~L.~R. Pereira, ``A
  three-phase optimal power-flow algorithm to mitigate voltage unbalance,''
  \emph{IEEE Trans. Power Del.}, vol.~28, no.~4, pp. 2394--2402, Oct. 2013.

\bibitem{DallaneseZhuGiannakis2013}
E.~Dall'Anese, H.~Zhu, and G.~B. Giannakis, ``Distributed optimal power flow
  for smart microgrids,'' \emph{IEEE Trans. Smart Grid}, vol.~4, no.~3, pp.
  1464--1475, Sept. 2013.

\bibitem{ZamzamSidiropoulosDallanese2016}
A.~S. Zamzam, N.~D. Sidiropoulos, and E.~Dall'Anese, ``Beyond relaxation and
  {Newton-Raphson}: Solving {AC} {OPF} for multi-phase systems with
  renewables,'' \emph{IEEE Trans. Smart Grid}, vol.~PP, no.~99, 2016.

\bibitem{WuKumagai1982}
F.~Wu and S.~Kumagai, ``Steady-state security regions of power systems,''
  \emph{IEEE Trans. Circuits Syst.}, vol.~29, no.~11, pp. 703--711, Nov. 1982.

\bibitem{Overbye1994}
T.~J. Overbye, ``A power flow measure for unsolvable cases,'' \emph{IEEE Trans.
  Power Syst.}, vol.~9, no.~3, pp. 1359--1365, Aug. 1994.

\bibitem{bolognani2016}
S.~Bolognani and S.~Zampieri, ``On the existence and linear approximation of
  the power flow solution in power distribution networks,'' \emph{IEEE Trans.
  Power Syst.}, vol.~31, no.~1, pp. 163--172, Jan. 2016.

\bibitem{YuTuritsyn2015}
S.~Yu, H.~D. Nguyen, and K.~S. Turitsyn, ``Simple certificate of solvability of
  power flow equations for distribution systems,'' in \emph{Proc. IEEE PES
  General Meeting}, July 2015, pp. 1--5.

\bibitem{WangBernsteinBoudecPaolone2016}
\BIBentryALTinterwordspacing
C.~Wang, A.~Bernstein, J.~Y.~L. Boudec, and M.~Paolone, ``Explicit conditions
  on existence and uniqueness of load-flow solutions in distribution
  networks,'' \emph{IEEE Trans. Smart Grid}, vol.~PP, no.~99, 2016. [Online].
  Available: \url{http://arxiv.org/abs/1602.08372}
\BIBentrySTDinterwordspacing

\bibitem{WangBernsteinBoudecPaolone2016threephase}
------, ``Existence and uniqueness of load-flow solutions in three-phase
  distribution networks,'' \emph{IEEE Trans. Power Syst.}, vol.~PP, no.~99,
  2016.

\bibitem{bazrafshanGatsis2016}
\BIBentryALTinterwordspacing
M.~Bazrafshan and N.~Gatsis, ``Convergence of the {Z-Bus} method for
  three-phase distribution load-flow with {ZIP} loads,'' \emph{IEEE Trans.
  Power Syst.}, May 2016, to be published. [Online]. Available:
  \url{https://arxiv.org/abs/1605.08511}
\BIBentrySTDinterwordspacing

\bibitem{RobbinsZhuGarcia2016}
B.~A. Robbins, H.~Zhu, and A.~D. Dominguez-Garcia, ``Optimal tap setting of
  voltage regulation transformers in unbalanced distribution systems,''
  \emph{IEEE Trans. Power Syst.}, vol.~31, no.~1, pp. 256--267, Jan. 2016.

\bibitem{BaranFernandes2016}
A.~R. {Baran Jr.} and T.~S.~P. Fernandes, ``A three-phase optimal power flow
  applied to the planning of unbalanced distribution networks,'' \emph{Int. J.
  of Electrical Power \& Energy Systems}, vol.~74, pp. 301--309, 2016.

\bibitem{JiangBaldick1996}
D.~Jiang and R.~Baldick, ``Optimal electric distribution system switch
  reconfiguration and capacitor control,'' \emph{IEEE Trans Power Syst.},
  vol.~11, no.~2, pp. 890--897, May 1996.

\bibitem{BolognaniDorfler2015}
S.~Bolognani and F.~D\"{o}rfler, ``Fast power system analysis via implicit
  linearization of the power flow manifold,'' in \emph{53rd Annu. Allerton
  Conf. Communication, Control, and Computing}, Sept. 2015, pp. 402--409.

\bibitem{KekatosZhangGiannakisBaldick2016}
V.~Kekatos, L.~Zhang, G.~B. Giannakis, and R.~Baldick, ``Voltage regulation
  algorithms for multiphase power distribution grids,'' \emph{IEEE Trans. Power
  Syst.}, vol.~31, no.~5, pp. 3913--3923, Sept 2016.

\bibitem{AhmadiMartiMeier2016}
H.~Ahmadi, J.~R. Martı´, and A.~von Meier, ``A linear power flow formulation
  for three-phase distribution systems,'' \emph{IEEE Trans. Power Syst.},
  vol.~31, no.~6, pp. 5012--5021, Nov 2016.

\bibitem{Garces2016}
A.~Garces, ``A linear three-phase load flow for power distribution systems,''
  \emph{IEEE Trans. Power Syst.}, vol.~31, pp. 827--828, Jan. 2016.

\bibitem{Gorman1992}
M.~J. Gorman and J.~J. Grainger, ``Transformer modelling for distribution
  system studies part {II}: Addition of models to {YBUS} and {ZBUS},''
  \emph{IEEE Trans. Power Del.}, vol.~7, no.~2, pp. 575--580, Apr. 1992.

\bibitem{chendillon1974}
M.-S. Chen and W.~E. Dillon, ``Power system modeling,'' \emph{Proc. of the
  {IEEE}}, vol.~62, no.~7, pp. 901--915, July 1974.

\bibitem{KerstingBook2001}
W.~H. Kersting, \emph{Distribution System Modeling and Analysis}, 3rd~ed.\hskip
  1em plus 0.5em minus 0.4em\relax CRC Press, 2002.

\bibitem{Moorthy2002}
S.~S. Moorthy and D.~Hoadley, ``A new phase-coordinate transformer model for
  {Ybus} analysis,'' \emph{IEEE Trans. Power Syst.}, vol.~17, no.~4, pp.
  951--956, Nov. 2002.

\bibitem{Kersting1999}
W.~H. Kersting, W.~H. Philips, and W.~Carr, ``A new approach to modeling
  three-phase transformer connections,'' \emph{IEEE Trans. Ind. Appl.},
  vol.~35, no.~1, pp. 169--175, Jan. 1999.

\bibitem{Xiao2006}
P.~Xiao, D.~C. Yu, and W.~Yan, ``A unified three-phase transformer model for
  distribution load flow calculations,'' \emph{IEEE Trans. Power Syst.},
  vol.~21, no.~1, pp. 153--159, Feb. 2006.

\bibitem{Dzafic2015}
I.~D\u{z}afi\'{c}, R.~A. Jabr, and H.-T. Neisius, ``Transformer modeling for
  three-phase distribution network analysis,'' \emph{IEEE Trans. Power Syst.},
  vol.~30, no.~5, pp. 2604--2611, Sept. 2015.

\bibitem{AndersonWollenberg1995}
D.~M. Anderson and B.~F. Wollenberg, ``Solving for three phase conductively
  isolated busbar voltages using phase component analysis,'' \emph{IEEE Trans.
  Power Syst.}, vol.~10, no.~1, pp. 98--108, Feb. 1995.

\bibitem{Chen1991}
T.~H. Chen, M.~S. Chen, T.~Inoue, P.~Kotas, and E.~A. Chebli, ``Three-phase
  cogenerator and transformer models for distribution system analysis,''
  \emph{IEEE Trans. Power Del.}, vol.~6, no.~4, pp. 1671--1681, Oct. 1991.

\bibitem{ieeefeederdata}
\BIBentryALTinterwordspacing
{IEEE PES Distribution System Analysis Subcommittee's Distribution Test Feeder
  Working Group}, ``Distribution test feeders.'' [Online]. Available:
  \url{http://ewh.ieee.org/soc/pes/dsacom/testfeeders/index.html}
\BIBentrySTDinterwordspacing

\bibitem{loadmodels1995}
``Standard load models for power flow and dynamic performance simulation,''
  \emph{{IEEE} Trans. Power Syst.}, vol.~10, no.~3, pp. 1302--1313, 1995.

\bibitem{loadmodels1995biblio}
``Bibliography on load models for power flow and dynamic performance
  simulation,'' \emph{{IEEE} Trans. Power Syst.}, vol.~10, no.~1, pp. 523--538,
  1995.

\bibitem{McKenna_2016}
K.~McKenna and A.~Keane, ``Open and closed-loop residential load models for
  assessment of conservation voltage reduction,'' \emph{{IEEE} Trans. Power
  Syst.}, pp. 1--1, 2016.

\bibitem{Chen1992}
T.~H. Chen and J.~D. Chang, ``Open wye-open delta and open delta-open delta
  transformer models for rigorous distribution system analysis,'' \emph{Proc.
  Inst. Elect. Eng. C - Generation Transmission Distribution}, vol. 139, no.~3,
  pp. 227--234, May 1992.

\bibitem{sallam2011electric}
A.~A. Sallam and O.~P. Malik, \emph{Electric distribution systems}.\hskip 1em
  plus 0.5em minus 0.4em\relax John Wiley \& Sons, 2011, vol.~68.

\bibitem{openDSSManual}
R.~C. Dugan, \emph{{OpenDSS Manual}}, EPRI, March 2016.

\bibitem{openDSS}
\BIBentryALTinterwordspacing
{EPRI}, ``{The Open Distribution System Simulator\textsuperscript{TM}, OpenDSS
  Version 7.6.5}.'' [Online]. Available:
  \url{https://sourceforge.net/projects/electricdss/}
\BIBentrySTDinterwordspacing

\bibitem{technoteTransformers}
{OpenDSS Wiki}, ``Technote transformers,'' OpenDSS/Doc/TechNotes.

\bibitem{IEEE8500}
R.~C. Dugan and R.~F. Arritt, ``The {IEEE} 8500-node test feeder,''
  https://ewh.ieee.org/soc/pes/dsacom/testfeeders/, EPRI, April 2010.

\end{thebibliography}
}

\vspace*{\baselineskip}

\end{document}